\documentclass[12pt,leqno]{amsart}
\usepackage{amsmath,stmaryrd}
\usepackage{amssymb}
\usepackage{amsthm}
\usepackage{epsf}
\usepackage{graphicx}
\usepackage{esint} 
\usepackage{dsfont}
\usepackage[pagebackref]{hyperref} 
\hypersetup{pdfpagemode=FullScreen,  colorlinks=true} 
\usepackage{enumerate} 
\usepackage{xcolor,bbm,easybmat,bm}

\numberwithin{equation}{section}

\newcommand\footnoteref[1]{\protected@xdef\@thefnmark{\ref{#1}}\@footnotemark}
\makeatother

\newcommand{\vp}{\varphi}
\newcommand{\dr}{\partial}

\DeclareMathOperator{\Jac}{Jac}

\DeclareMathOperator{\diver}{div}

\DeclareMathOperator{\dist}{dist}
\DeclareMathOperator{\supp}{supp}
\DeclareMathOperator{\diam}{diam}

\def\div{\mathop{\operatorname{div}}}
\DeclareMathOperator{\loc}{loc}

\newcommand{\1}{{\mathds 1}}

\newcommand{\ms}{\medskip}

\newcommand{\R}{\mathbb R}
\newcommand{\N}{\mathbb N}

\newcommand{\bp}{\noindent {\em Proof: }}
\newcommand{\ep}{\hfill $\square$ \medskip}

\newcommand{\wt}{\widetilde}

\newcommand{\D}{\mathbb D}

\newcommand{\A}{\mathcal A}
\newcommand{\C}{\mathcal C}
\newcommand{\G}{\mathcal G}
\newcommand{\B}{\mathcal B}
\newcommand{\cS}{\mathcal S}
\newcommand{\W}{\mathcal W}

\newcommand{\Rn}{\mathbb R^n}
\newcommand{\norm}[1]{\left\Vert#1\right\Vert}
\newcommand{\abs}[1]{\left\vert#1\right\vert}
\newcommand{\br}[1]{\left(#1\right)}
\newcommand{\set}[1]{\left\{#1\right\}}

\newcommand{\om}{\Omega}
\newcommand{\pom}{\partial\Omega}

\newcommand{\dint}{\int\!\!\!\!\!\int}
\def\Yint#1{\mathchoice
	{\YYint\displaystyle\textstyle{#1}}%
	{\YYint\textstyle\scriptstyle{#1}}%
	{\YYint\scriptstyle\scriptscriptstyle{#1}}%
	{\YYint\scriptscriptstyle\scriptscriptstyle{#1}}%
	\!\dint}
\def\YYint#1#2#3{{\setbox0=\hbox{$#1{#2#3}{\iint}$}
		\vcenter{\hbox{$#2#3$}}\kern-.51\wd0}}
\def\longdash{\mkern-1.5mu{-}\mkern-7.5mu{-}} 

\def\fiint{\Yint\longdash}

\theoremstyle{plain}
\newtheorem{theorem}[equation]{Theorem}
\newtheorem{lemma}[equation]{Lemma}
\newtheorem{corollary}[equation]{Corollary}
\newtheorem{proposition}[equation]{Proposition}

\newtheorem{definition}[equation]{Definition}

\theoremstyle{definition}

\theoremstyle{remark}
\newtheorem{remark}[equation]{Remark}

\begin{document}

\title[A Green function characterization of uniformly rectifiable sets]{A Green function characterization of uniformly rectifiable sets of any codimension}

\author[Feneuil]{Joseph Feneuil}
\address{Joseph Feneuil. Mathematical Sciences Institute, Australian National University, Acton, ACT, Australia}
\email{joseph.feneuil@anu.edu.au}

\author[Li]{Linhan Li}
\address{Linhan Li. School of Mathematics, The University of Edinburgh, Edinburgh, UK }
\email{linhan.li@ed.ac.uk}

\thanks{
Early discussions on this project took place during the first author's stay in Minneapolis - funded by the Simons Foundation grant 563916, SM - 
and the major part of this work was carried out during the period when the second author was visiting the Australian National University. The authors would like to express their gratitude to Svitlana Mayboroda, Po Lam Yung, and the Australian National University for their support.}

\maketitle

\begin{abstract} 
In this paper, we obtain a unified characterization of uniformly rectifiable sets of {\it any codimension} in terms of a Carleson estimate on the second derivatives of the Green function. When restricted to domains with boundaries of codimension 1, our result generalizes a previous result of Azzam for the Laplacian to more general elliptic operators. For domains with boundaries of codimension greater than 1, our result is completely new. 
\end{abstract}

\ms\noindent{\bf Keywords: uniform rectifiability, Green function, degenerate operators}

\ms\noindent

\tableofcontents

\section{Introduction}
\label{S1}

\subsection{State of the art}

Rectifiable sets play an important role in geometric measure theory 
and the calculus of variation, in particular
because the sets that minimize an energy are often rectifiable. In the past thirty decades, the scale-invariant and quantitative version of rectifiability- uniform rectifiability- has attracted a great deal of interest for its rich interactions with areas in and beyond geometric measure theory. The notion of uniform rectifiability was introduced by David and Semmes (\cite{DS1,DS2}) in the early 90's, 
along with many characterizations in terms of geometry and in terms of boundedness of singular integral operators. In recent years, a finale of a vast amount of work brought a characterization of uniform rectifiability in terms of elliptic PDEs. More precisely, it has been shown that if a domain $\Omega$ is sufficiently well connected, the uniform rectifiability of $\partial \Omega$ can be characterized by (quantitatively) absolute continuity of the harmonic measure with respect to the surface measure, by the $L^p$ solvability of the Dirichlet problem for a large $p$, as well as by Carleson estimates for solutions. A lot of research and many mathematicians were needed to reach these PDE characterizations, a non-exhaustive list of works includes \cite{DJ90,HM1,HMUT14,AHMNT17,AzzamTravelingSaleman}. The optimal topological conditions are obtained in \cite{Azzam21} (for the $A_\infty$-absolute continuity of the harmonic measure) and \cite{AHM3T20} (for the solvability of the Dirichlet problem in $L^p$).
There have been also beautiful and difficult extensions of these characterizations to more general elliptic operators, for example, \cite{HMT,HMMtrans,HMMTZ}, to just mention a few.

Nevertheless, the above mentioned PDE characterizations of uniform rectifiability are limited to sets that are seen by the Laplacian, that is, the boundary of the underlying domain has to have codimension one, while uniform rectifiability is a notion that exists for any integer dimension. 
In \cite{DFMhighcod,DFMmixed}, David, Mayboroda, and the first author developed an elliptic theory in the complement of a set $E$ of codimension greater than one with the help of degenerate elliptic operators. It has been proved that the elliptic measure on $E$ of a carefully chosen degenerate operator is $A_\infty$-absolutely continuous with respect to the Hausdorff measure on $E$ when $E$ is  uniformly rectifiable (see \cite{DFMdahlberg, DMAinfty, FenAinfty}).
However, it was discovered in \cite[Theorem 6.7]{DEM} that the converse is false for one of these operators. 
This failure of characterizing uniformly rectifiable sets of higher codimension by absolute continuity of harmonic measure leads to intensive investigations on other PDE properties, in particular, properties of the Green function (see \cite{DMgreen,DFM22,DLM1,DLMlower,FLM}). Some beautiful and groundbreaking results have been obtained and some of them are very different in form and nature, which
we shall discuss in more detail later in Section~\ref{subsecComp}, but a unified characterization of uniform rectifiability in all dimensions in terms of a strong Carleson condition has yet to be established.    

In this paper, building on earlier works (mainly \cite{DFM22,FLM,DEM}), we obtain a desired characterization of uniformly rectifiable sets of {\it any codimension} in terms of a Carleson estimate on the second derivatives of the Green function. Our characterization can take the same form in all codimensions and can be viewed as an extension of \cite{AzzamTravelingSaleman} from the Laplacian to more general elliptic operators, and to domains with a boundary of codimension greater than one. 
\smallskip

Let us be more precise. The statement of \cite[Theorem VI]{AzzamTravelingSaleman} entails
\footnote{The formulation in \cite{AzzamTravelingSaleman} is slightly different, and the case where the boundary is lower $(n-1)$-content regular, a condition weaker than  ($n-1$)-Ahlfors regular, is also considered. But to get a characterization of uniform rectifiability one has to assume ($n-1$)-Ahlfors regularity.}
that for a bounded uniform domain $\om$ in $\Rn$ with ($n-1$)-Ahlfors regular boundary, $\pom$ is uniformly rectifiable (see Definition \ref{defUR}) if and only if the Green function $G^Y$ for the Laplacian with pole at some point $Y\in\om$ satisfies the estimate 
\[
\iint_{B(x,r)\cap\om}\abs{\frac{\nabla^2G^Y(X)}{G^Y(X)}}^2\dist(X,\pom)^3dX\le Cr^{n-1}
\]
for any $x\in\pom$ and $r\in(0,\diam\pom)$ such that $Y\notin B(x,2r)$. 
Observe that the quantity $\abs{\frac{\nabla^2G^Y(X)}{G^Y(X)}}\dist(X,\pom)^2$
measures how close the Green function $G$ is to being affine around $X$ in a scale-invariant way. Although this result is formulated for the Laplacian only, using a similar technique in \cite{HMT}, one can probably extend it to a slightly more general class of elliptic operators. Nonetheless, with this method, the optimal class of elliptic operators (see Definition \ref{def.DKPcd}), which was introduced by Dahlberg and proved by Kenig and Pipher in \cite{KP01} to be the right class of elliptic operators to consider for absolute continuity of elliptic measure in Lipschitz domains, is still far from reach. Our result, when restricted to the case of codimension 1, applies to this optimal class of elliptic operators. In the case where the boundary of the domain is of codimension greater than 1, our result is completely new.

\subsection{The statement of the main results}

We take a domain $\Omega \subset \R^n$ whose boundary $\partial \Omega$ is $d$-Ahlfors regular with $0<d\le n-1$, which means that there exists a measure $\sigma$ supported on $\partial \Omega$ and a constant $C=C_\sigma>0$ such that 
\begin{equation} \label{defADR}
C_\sigma^{-1}r^d \leq \sigma(B(x,r)) \leq C_\sigma r^d \qquad \text{ for } x\in \partial \Omega, \, r \in (0, \diam \partial \Omega).
\end{equation}
Note that we can always take $\sigma$ to be the $d$-dimensional Hausdorff measure on $\partial \Omega$. 
\medskip

\medskip

We also need a smooth distance function that was introduced in \cite{DFMdahlberg} and is defined, for some $\beta>0$, as
\begin{equation} \label{defDbeta}
D_\beta(X) := \left( \int_{\partial \Omega} |X-y|^{-d-\beta} d\sigma(y) \right)^{-\frac1\beta} \qquad \text{ for $X\in \Omega$}.
\end{equation}
One can show that when $\partial \Omega$ is $d$-Ahlfors regular (and $\beta >0$), the smooth distance $D_\beta$ is equivalent to $\dist(\cdot,\pom)$ - which we denote by $\delta_{\partial \Omega}$ from now on - that is, for some $C_\beta\ge1$,
\begin{equation} \label{Db=dist}
C_\beta^{-1} \delta_{\partial \Omega}(X) \leq D_\beta(X)  \leq C_\beta \delta_{\partial \Omega}(X)
\end{equation}
whenever $X\in \Omega$ satisfies $\delta_{\partial \Omega}(X) \leq 100 \diam(\partial \Omega)$, see \cite[Lemma 5.1]{DFMdahlberg}.
\medskip

We say that $L:= - \diver A \nabla$ is a uniformly elliptic\footnote{Note that when $d<n-1$, this is not the standard definition and is usually referred to as degenerate elliptic. We avoid this complication in this paper.} operator
if there exists $C>0$ - called ``{\bf elliptic constant of $L$}'' in our paper - such that 
\begin{equation} \label{defelliptic2}
A(X)\xi\cdot \xi \geq C^{-1} \delta_{\partial \Omega}(X)^{d+1-n} |\xi|^2 \qquad \text{ for } X\in \Omega, \, \xi \in \R^n,
\end{equation}
and
\begin{equation} \label{defbounded2}
|A(X)\xi \cdot \zeta| \leq C \delta_{\partial \Omega}(X)^{d+1-n} |\xi| |\zeta| \qquad \text{ for } X\in \Omega, \, \xi \in \R^n.
\end{equation}
In virtue of \eqref{Db=dist}, the weight is equivalent to $D_\beta(X)^{d+1-n}$, so in the paper, we say that $L=-\diver(D_\beta^{d+1-n} \A \nabla)$ is uniformly elliptic if $\A$ satisfies the classical elliptic and boundedness conditions:
\begin{equation} \label{defelliptic}
\A(X)\xi\cdot \xi \geq C^{-1} |\xi|^2 \qquad \text{ for } X\in \Omega, \, \xi \in \R^n,
\end{equation}
and
\begin{equation} \label{defbounded}
|\A(X)\xi \cdot \zeta| \leq C |\xi| |\zeta| \qquad \text{ for } X\in \Omega, \, \xi \in \R^n.
\end{equation}
Let $B \subset \R^n$ be a ball centered on the boundary or inside $\Omega$. A (weak) solution to $Lu = -\div A \nabla u =0$ in $B\cap \Omega$ means that $u$ belongs to 
\[W(B\cap \Omega) := \Big\{v \in W^{1,2}_{loc}(B\cap \Omega), \, \iint_{B \cap \Omega} |\nabla v|^2 \, \delta_{\partial \Omega}^{d+1-n}\, dX<\infty \Big\}\]
and satisfies
\[\iint_{\Omega} A \nabla u \cdot \nabla \varphi \, dX = 0 \qquad \text{ for } \varphi \in C^\infty_0(B\cap \Omega).\] 
If, in addition, the domain $\Omega$ is uniform (see Definition \ref{defUniform}) - which is always the case in our paper - then $u$ is continuous on $B\cap \overline{\Omega}$ (see the elliptic theory developed in \cite{DFMhighcod}--\cite{DFMmixed}, in particular, \cite[Lemmas 11.30 and 11.32]{DFMmixed}), and the expression $u = 0$ on $\partial \Omega \cap B$ makes sense.

\bigskip

The notion of Carleson measure
will be central to our paper. 
\begin{definition}
We say that a quantity $f$ defined on $\Omega$ satisfies the Carleson measure condition - or $f\in CM_\Omega$ for short - if there exists a constant $M>0$ such that for any $x\in \dr \Omega$ and $r<\diam(100 \partial \Omega)$,
\begin{equation} \label{defCarleson}
\iint_{B(x,r) \cap \Omega} |f(X)|^2 \delta(X)^{d-n} dX \leq M^2 r^{d}.
\end{equation}
When we want to emphasize the constant $M$, we also write $f\in CM_\om(M)$.
\end{definition}

\begin{definition}[The Dahlberg-Kenig-Pipher (DKP) condition]\label{def.DKPcd}
Let $\A$ be a matrix-valued function defined on $\om$ that satisfies \eqref{defelliptic} and \eqref{defbounded}, and let $L=-\diver\A\nabla$. We say that the matrix $\A$ or the operator $L$ satisfies the Dahlberg-Kenig-Pipher (DKP) condition if there exists $M>0$, such that 
\begin{equation}\label{DAbdd}
        D(X)\abs{\nabla\A(X)}\le M \quad \text{for all }X\in\om,
\end{equation}
and 
\begin{equation} \label{defDKP}
     D|\nabla\A| \in CM_\Omega(M),
\end{equation}
where $D$ stands for either $\delta_{\partial \Omega}$ or $D_\beta$ for some $\beta>0$. We call the constant $M$ the `` {\bf DKP constant of $L$}''. If the operator $L$ satisfies the DKP condition, we also say it is a DKP-operator. 
\end{definition}

Note that the DKP condition will be considered only in the case of codimension 1 boundaries ($d=n-1$) in this paper. However, similar conditions have been considered in higher codimension cases (e.g. \cite{DMgreen,DLMlower}). In general, the analogous condition will be more complicated in the higher codimensional settings. 
 Our proof should include this case provided that one of the results that we relied on, namely, the main result in \cite{DFM22}, is established for DKP-operators in higher codimension.

\medskip

If we set aside the technical definition of uniform domain and uniform rectifiability for now (see Definitions \ref{defUniform} and \ref{defUR} for the precise statements), we are ready to present our results. We start with the codimension one case.

\begin{theorem}[Domains with codimension 1 boundaries] \label{Maincod1}
Let $\Omega \subset \R^n$ be a uniform domain with $(n-1)$-Ahlfors regular boundary, and let $L = - \div \A \nabla$ be a DKP-operator. Then the following are equivalent:
\begin{enumerate}[(i)]
\item $\partial \Omega$ is uniformly rectifiable;
\item there exists $C>0$ such that for any $x\in \dr \Omega$, any $r<\diam(100 \partial \Omega)$, and any positive weak solution $u$ to $Lu = 0$ in $\Omega \cap B(x,2r)$ satisfying $u = 0$ in $\partial \Omega \cap B(x,2r)$, we have
\begin{equation} \label{N2u/uisCM}
\iint_{B(x,r)\cap\om} \frac{|\nabla^2 u(X)|^2}{u(X)^2} \delta_{\partial \Omega}(X)^{3} dX \leq C r^{d};
\end{equation}
\item there exists $C>0$ and $Y\in \Omega$\footnotemark{}
\footnotetext{we also assume $\delta_{\partial \Omega}(Y) \leq 100 \diam\partial \Omega$ if $\diam \partial \Omega < \diam \Omega = \infty$.}
such that for any $x\in \dr \Omega$ and any $r>0$ satisfying $Y \notin B(x,2r)$, we have
\begin{equation} \label{D2G/GCM}
\iint_{B(x,r)} \frac{\abs{\nabla^2 G^Y(X)}^2}{G^Y(X)^2} \delta_{\partial \Omega}(X)^{3} dX \leq C r^{d};
\end{equation}
\item there exists $C>0$ and $Y\in \Omega$\footnotemark[\value{footnote}]
such that for any $x\in \dr \Omega$ and any $r>0$ satisfying $Y \notin B(x,2r)$, we have
\begin{equation} \label{NNu/uisCM}
\iint_{B(x,r)} \frac{\abs{\nabla \br{|\nabla G^Y(X)|^2}}^2}{G^Y(X)^4} \delta_{\partial \Omega}(X)^{5} dX \leq C r^{d},
\end{equation}
where $G^Y$ is the Green function of $G$ with pole at $Y$ (see Definition \ref{defGreen}).
\end{enumerate}
\end{theorem}

We shall prove $(i)\implies(ii)\implies(iii)\implies(iv)\implies(i)$. Note that $(ii)\implies(iii)$ is immediate as $G^Y$ is a positive solution in $B(x,2r)\cap\om$ (when $Y\notin B(x,2r)$) and vanishes on $B(x,2r)\cap\pom$. The implication $(iii)\implies(iv)$ is also straightforward because 
\[
    \frac{\delta_{\pom}\abs{\partial_i \br{|\nabla G|^2}}}{G}= \frac{\delta_{\pom}\abs{\partial_i \br{\nabla G\cdot\nabla G}}}{G}\le \frac{2\delta_{\pom}\abs{\nabla G}}{G}\abs{\nabla^2G}\le C\abs{\nabla^2G}
\]
by Lemma \ref{lemDNu/u<1} below. 

\medskip

As mentioned earlier, we aim to obtain a unified characterization of uniformly rectifiable sets of all codimension. We let the reader check that \eqref{N2u/uisCM} and \eqref{D2G/GCM} have to be specific to codimension 1 boundary, since taking the simple case $\Omega = \R^n\setminus\R^d$, $L= - \div |t|^{d+1-n} \nabla$, and $u = |t|$ gives a counterexample. It turns out that expressions like \eqref{NNu/uisCM} work for domains with boundaries of any codimension.

\begin{theorem}
\label{mainthm}
  Let $0<d<n$ and $\Omega$ be a domain in $\R^n$ with $d$-Ahlfors regular boundary. Take $\beta >0$ and define $D_\beta$ as in \eqref{defDbeta}. Let $L=-\diver\br{D_\beta^{d+1-n}\A\nabla}$ be uniformly elliptic. Assume either
\begin{enumerate}
    \item $d=n-1$, $\Omega$ is uniform, and $\A$ satisfies the DKP condition, or 
    \item $0<d<n-1$ and $\A=I$.
\end{enumerate}
Then the following are equivalent:
\begin{enumerate}[(i)]
\item $d$ is an integer and $\partial \Omega$ is uniformly rectifiable;
\item there exists $C>0$ such that for any $x\in \dr \Omega$, any $r<\diam(100 \partial \Omega)$, and any positive weak solution $u$ to $Lu = 0$ in $\Omega \cap B(x,2r)$ satisfying $u = 0$ in $\partial \Omega \cap B(x,2r)$, we have
\begin{equation} \label{N2u/uisCMb}
\iint_{B(x,r)} \frac{\big|\nabla (|\nabla u(X)|)\big|^2}{u(X)^2} \delta_{\partial \Omega}(X)^{d+4-n} dX \leq C r^{d}.
\end{equation}
\item there exists $C>0$ such that for any $x\in \dr \Omega$, any $r<\diam(100 \partial \Omega)$, and any positive weak solution $u$ to $Lu = 0$ in $\Omega \cap B(x,2r)$ satisfying $u = 0$ in $\partial \Omega \cap B(x,2r)$, we have
\begin{equation} \label{N2u/uisCMc}
\iint_{B(x,r)} \frac{\big|\nabla \br{|\nabla u(X)|^2}\big|^2}{u(X)^4} \delta_{\partial \Omega}(X)^{d+6-n} dX \leq C r^{d};
\end{equation}
\item there exists $C>0$ and $Y\in \Omega$\footnotemark[\value{footnote}]
such that for any $x\in \dr \Omega$ and any $r>0$ satisfying $Y \notin B(x,2r)$, we have
\begin{equation} \label{NNu/uisCMb}
\iint_{B(x,r)} \frac{\abs{\nabla \br{|\nabla G^Y(X)|^2}}^2}{G^Y(X)^4} \delta_{\partial \Omega}(X)^{d+6-n} dX \leq C r^{d},
\end{equation}
where $G^Y$ is the Green function of $G$ with pole at $Y$.
\end{enumerate}
\end{theorem}

We shall show $(i)\implies(ii)\implies(iii)\implies(iv)\implies(i)$. The implication $(iii)\implies(iv)$ is trivial, and $(ii)\implies(iii)$ follows from
\[ \frac{\delta_{\partial \Omega} \big|\nabla |\nabla u|^2 \big|}{u} = \frac{2\delta_{\partial\Omega} |\nabla u|}{u} \big|\nabla |\nabla u| \big| \leq C\big|\nabla |\nabla u| \big|,\]
where in the last inequality we have used again Lemma \ref{lemDNu/u<1}.

\begin{remark}
Note that our method to show the free boundary argument $(iv)\implies(i)$ can be adapted to consider a larger class of operators. In Theorem \ref{Thconverse}, we will show that a weak version of $(iv)$ implies $(i)$ for operators in the form $L= - \diver \br{D_\beta^{d+1-n} \A\nabla}$, where $\A$ satisfies the DKP condition. We claim that we could also prove $(iv)\implies(i)$ for operators in the form $L= - \diver \br{D_\beta^{d+1-n} \A\nabla }$ where $\A \in C^\infty(\Omega)$ and $\delta_{\partial \Omega} \nabla \A \in L^\infty(\Omega)$. An example of non-DKP operator satisfying the latter condition is the one constructed in the complement of the 4 corners Cantor set by David and Mayboroda (see \cite{DM4corners}), hence we can obtain that $(ii)$--$(iv)$ are necessarily false for this operator from \cite{DM4corners}. 
\end{remark}

\medskip

\subsection{Comparison to earlier works and ingredients of the proof}\label{subsecComp}

As we mentioned earlier, on domains with boundaries of codimension 1, formulations similar to our Theorem \ref{Maincod1} have shown up in different contexts. The reader should have recognized that \cite[Theorem VI]{AzzamTravelingSaleman} proves Theorem \ref{Maincod1} $(i)\Longleftrightarrow(iii)$ when $L$ is the Laplacian. The proof of Propositions 4.25 and 4.36 in \cite{HMT} contains the arguments that would prove $(i) \implies (iii)$ when $L=-\div \A \nabla$ satisfies \eqref{DAbdd} and $(D|\nabla \A|)^{1/2} \in CM_\Omega$, which is stronger than our condition $D|\nabla \A| \in CM_\Omega$. 
 The proof of the results in \cite{HMT} and \cite{AzzamTravelingSaleman} use a similar argument that could not be easily adapted to include the full class of DKP operators.
By a completely different method, \cite[Corollary 1.17]{DLM1} proves $(i) \implies (ii)$ for the full class of DKP operators, but only when $\Omega$ is the half space (it immediately extends to Lipschitz domains, but we do not see how the method of \cite{DLM1} can be used for domains with uniformly rectifiable boundaries).
\medskip

Our proof uses a third strategy, which utilizes a smooth distance $D_\beta$ that does not show up in the statements of our theorems at all. In \cite{DEM}, David, Engelstein, and Mayboroda showed that the oscillation of $\abs{\nabla D_\beta}$ characterizes  uniformly rectifiable sets of any codimension. To be precise, they proved the following result.

\begin{theorem}[\cite{DEM} Theorem 1.4]\label{thmDEM}
Let $n\ge 2$ and $0<d<n$. A $d$-Ahlfors regular
set $E$ in $\R^n$ equipped with a $d$-Ahlfors regular measure $\sigma$ is uniformly rectifiable if and
only if there exists $M>0$ and $\beta>0$ such that 
\begin{equation} \label{dem est}
 \delta_E\abs{\nabla\br{\abs{\nabla D_{\beta}}^2}}\in CM_{\Rn\setminus E}(M).
\end{equation}
The result still holds when replacing \eqref{dem est} by
\begin{equation} \label{dem est 2}
 \delta_E\abs{\nabla\br{\abs{\nabla D_{\beta}}}}\in CM_{\Rn\setminus E}(M).
\end{equation}
\end{theorem}
We shall only need the ``direct" direction of this result, that is, if $E$ is uniformly rectifiable, then \eqref{dem est} and \eqref{dem est 2} hold. However, to prove $(iv)\implies(i)$ in Theorem \ref{mainthm}, we borrow ideas from the proof of the other direction of Theorem \ref{thmDEM}. 

Another main ingredient in our proof is a comparison between the Green function and the smooth distance $D_\beta$. In the recent work \cite{FLM}, the authors and Mayboroda gave a characterization of uniformly rectifiable sets of codimension 1 by a Carleson estimate on $\nabla\ln\br{G/D_\beta}$. More precisely, we showed that if $\om\subset\Rn$ is a uniform domain with $(n-1)$-Ahlfors regular boundary, and if $L=-\div \A \nabla$ satisfies the (weak) DKP condition, then 
 $\partial \Omega$ is uniformly rectifiable if and only if there exist $Y\in \Omega$ and $C>0$  such that for any ball $B$ centered on the boundary that satisfies $Y\notin 2B$, there holds
    \begin{equation} \label{eqFLM}
    \iint_{\Omega \cap B} \left| \frac{\nabla G^{Y}}{G^{Y}} - \frac{\nabla D_\beta}{D_\beta} \right|^2 D_\beta \, dX \leq C r^{n-1}.
    \end{equation}
The idea of relating uniform rectifiability to a comparison between the Green function and a distance function emerged in works before \cite{FLM}, and in particular, \cite{DMgreen} and \cite{DFM22} consider the higher codimension case as well. Unfortunately, in the case when the boundary of the underlying domain has codimension strictly greater than 2, \eqref{eqFLM} fails to characterize uniform rectifiability as it has been observed in \cite{DEM} that when $n-d-2>0$, the smooth distance $D_\beta=D_{n-d-2}$ is a solution to $\diver( D_\beta^{d+1-n} \nabla u)=0$ regardless of the geometry of the boundary. So $D_{n-d-2}$ is the Green function ``with pole at infinity'', and thus \eqref{eqFLM} holds trivially (if taking $Y$ to be infinity) yet it does not tell anything about the boundary. This issue leads to two consequences in previous works. First, the characterization given in \cite{DMgreen} using the difference between the Green function and $D_\beta$ has to rule out the codimension 2 case and to assume a condition on the domain (called the property $\Upsilon_{\text{flat}}(d,\beta,\Delta)$ in \cite{DMgreen}) that people do not yet know when is valid (see \cite{Pe} for some partial progress). Second, the result in \cite{DFM22}, which deals with the higher codimension case only, holds only in one direction. Nonetheless, the comparison between the Green function and $D_\beta$ plays an important role in our proof. We record the results obtained in \cite{DFM22} and \cite{FLM} in the following theorem.

\begin{theorem}[\cite{FLM} Theorem 1.12, \cite{DFM22}  Theorem 2.21] 
\label{ThNlnu/D} 
Let $\Omega \subset \Rn$ be a domain with uniformly $d$-rectifiable boundary, $\beta >0$, and $L=-\diver\br{D_\beta^{d+1-n}\A\nabla}$. Assume either
\begin{enumerate}
    \item $d=n-1$, $\Omega$ is uniform, and $\A$ satisfies the DKP condition, or 
    \item $0<d<n-1$ and $\A=I$.
\end{enumerate}
Then for any $x\in \partial \Omega$, any $r\in (0,100\diam \partial \Omega)$, and any positive weak solution $u$ to $Lu = 0$ in $\Omega \cap B(x,2r)$ with $u=0$ on $\partial \Omega \cap B(x,2r)$, one has
\begin{equation} \label{GiCM}
\iint_{B(x,r)\cap\om} \left|\nabla \ln\left( \frac{u}{D_{\beta}} \right)\right|^2 \, D_\beta^{d+2-n} \,  dX \leq C r^{d}.
\end{equation}
\end{theorem}

To summarize, it is understood that when $\partial \Omega$ is uniformly rectifiable and the operator $L$ is ``suitable'', the Green function should behave  like $D_\beta$ when away from the pole. Then, we should be able to characterize the uniform rectifiability of boundaries of any codimension by replacing $D_\beta$ in Theorem \ref{thmDEM} by the Green function of ``suitable'' elliptic operators, which is our Theorem \ref{mainthm}.

\subsection{Organization of the article}
In Section~\ref{SecDef}, we give definitions of the geometric notions used in this paper and collect some regularity results for the solutions. 

Sections~\ref{SecD},\ref{SecCaccio},\ref{SecURtoG} are devoted to proving the ``direct" direction, that is, uniform rectifiability implies Carleson estimates for the Green function. 
We prove $(i) \implies (ii)$ in Theorems \ref{Maincod1} and \ref{mainthm} by first establishing a second derivative version of Theorem \ref{ThNlnu/D} (see Corollary \ref{CorCaccio}) in Section~\ref{SecCaccio} using a Caccioppoli type argument, which produces some terms involving the derivatives of $D_\beta$. Then in Section~\ref{SecURtoG}, we play with the second derivatives of the solution  and invoke Corollary \ref{CorCaccio}, estimates on $D_\beta$ and its derivatives, as well as Theorem \ref{ThNlnu/D} to get the desired result. 

The remaining sections are for the free boundary direction, namely, $(iv)\implies(i)$ in Theorem \ref{mainthm}. We use a blow-up argument in Section~\ref{SecGtoUR} to reduce the problem to the limit case, which is treated in Section~\ref{SecLimit}. Theorem \ref{mainthm} $(iv) \implies (i)$ is a consequence of Theorem \ref{Thconverse}, and the proof is given in the end of Section~\ref{SecGtoUR}.

 In the sequel, we shall use $A \lesssim B$ when $A \leq C B$ for a constant $C$ that depends only on the relevant parameters (that are either recalled or obvious from context). We write $A \approx B$ if $A \lesssim B$ and $B \lesssim A$.

\section{Definitions and existing results}\label{SecDef}
\subsection{Topological assumptions}

\begin{definition} \label{defCP}
A domain $\Omega \subset \R^n$ satisfies the {\bf corkscrew point condition} if there exists $\epsilon\in(0,1)$ such that for any $x\in \partial \Omega$, and any $r\in (0,\diam\Omega)$, there exists $X \in \Omega \cap B(x,r)$ such that $B(X,\epsilon r) \subset \Omega$.
\end{definition}

\begin{definition} \label{defHC}
Let $\Omega \subset \R^n$ be a domain that satisfies the corkscrew point condition (with constant $\epsilon$). We say that $\Omega$ satisfies the {\bf Harnack chain condition} if there exists $N$ such that for any couple of point $(X,Y) \in \Omega\times\om$ satisfying $|X-Y|/\min\{\delta_{\pom}(X),\delta_{\pom}(Y)\} \leq 10/\epsilon$, there exists $N+1$ points $\{Z_i\}_{0\leq i\leq N}$ such that $Z_0 = X$, $Z_{N} = Y$, and $|Z_i-Z_{i+1}| \leq \delta_{\pom}(Z_i)/2$.
\end{definition}

\begin{definition} \label{defUniform}
A domain $\Omega$ is said to be {\bf uniform}\footnote{a.k.a. 1-sided NTA domain} if $\Omega$ satisfies both the corkscrew point condition and the Harnack chain condition. The constants of a uniform domain $\Omega$ are the couple $(\epsilon,N)$ given in Definitions \ref{defCP} and \ref{defHC}.
\end{definition}

There are several equivalent definitions of uniform domains. We shall just mention here that a uniform domain is necessarily connected, and any couple of points can be linked by a thick path, as shown by the following proposition.

 \begin{proposition}\cite[Proposition 2.18]{DFMmixed} \label{propHarnack}
Let $\Omega$ be a uniform domain with constants $\epsilon,N$. 
There exists $C:=C(\epsilon,N)>0$ such that for any $\Lambda \geq 1$, and any couple $X,Y \in \Omega$ that satisfy $|X-Y|/\min\{\delta_{\pom}(X),\delta_{\pom}(Y)\} \leq \Lambda$, we can find $N_\Lambda:= \lceil C\ln(1+\Lambda)\rceil$ points $Z_0:=X,Z_1,\dots,Z_{N_\Lambda} = Y$ such that for any $i\in \{ 0,\dots,N_\Lambda-1\}$, 
\begin{enumerate}[(i)]
\item $|Z_i - Z_{i+1}| \leq \frac12\delta_{\pom}(Z_i)$, 
\item $\delta_{\pom}(Z_i) \geq 2^{-N}\min\{\delta_{\pom}(X),\delta_{\pom}(Y)\}$,
\item $|Z_i - X| + |Z_i - Y| \leq 2^{N+5}
\Lambda\min\set{\delta_{\pom}(X),\delta_{\pom}(Y)}$.
\end{enumerate}
\end{proposition}

Observe also that when $\Omega$ is the complement of a low dimensional set, then $\Omega$ is necessarily uniform.

\begin{proposition} \label{lowdim=>uniform}
Let $d<n-1$ and $\Omega \subset \R^n$ be a domain with $d$-Ahlfors regular boundaries. Then $\Omega$ is uniform.
\end{proposition}

\bp
It has been established in \cite[Lemmas 2.1 and 11.6]{DFMhighcod} that if we assume in addition that $\om=\Rn\setminus\pom$, then $\om$ is uniform. We claim that 
if $\partial \Omega$ is $d$-dimensional with $d<n-1$, then 
\begin{equation} \label{claimdO=R-O}
\Omega = \R^n \setminus \partial \Omega,
\end{equation} 
and this will complete the proof. 
Assume by contraposition that $\partial \Omega \subsetneq \R^n \setminus \Omega$, then we can find two small balls $B(x,r) \subset \R^n\setminus \Omega$ and $B(y,r) \subset \Omega$. By dilatation, translation, and rotation invariance, we can choose $x=0$, $y=(0,\dots,0,1)$, and $r<\frac12$. Denote by $\pi$ the orthogonal projection on $\R^{n-1} \times \{0\}$, which is $1$-Lipschitz. Observe that for any $z\in B_{\R^{n-1}}(0,r)$, the point $(z,1)\in B(y,r)$. So the line segment with end points $(z,0)$ and $(z,1)$ necessarily intersects $\pom$. Since $(z,t)\in B(0,2)$ for all $z\in B_{\R^{n-1}}(0,r)$ and $t\in[0,1]$, we have
 $B_{\R^{n-1}}(0,r)\subset \pi(\partial \Omega \cap B(0,2))$. As a consequence,
\[ \mathcal H^d(B_{\R^{n-1}}(0,r)) \leq \mathcal H^d(\pi(\partial \Omega \cap B(0,2))) \leq \mathcal H^d(\partial \Omega \cap B(0,2)) < \infty,\]
which necessarily means that $d\geq n-1$.
\ep

\subsection{Dyadic and Whitney decompositions} \label{SSdyadic} 

Let $\om$ be a domain in $\Rn$ with $d$-Ahlfors regular boundary.

\subsubsection{Dyadic decomposition}
We construct a dyadic system of pseudo-cubes on $\partial \Omega$. In the presence of the Ahlfors regularity property, such construction appeared for instance in \cite{David91}, \cite{DS1} or \cite{DS2}. As in \cite{FLM}, we shall use the construction of Christ \cite{Ch}, which allows to bypass the need of a measure on $\partial \Omega$. 

There exist a universal constant $0<a_0<1$ and a collection $\mathbb{D}_{\partial \Omega} = \bigcup_{k \in \mathbb Z} \mathbb{D}_{\partial \Omega, k}$ of Borel subsets of $\partial \Omega$ with the following properties. We write 
\[\mathbb{D}_{\partial \Omega, k}:=\{Q_{j}^k\subset \D_{\partial \Omega}: j\in \mathfrak{I}_k\}, \]
where $\mathfrak{I}_k$ denotes some ordered index set depending on $k$, but sometimes, to lighten the notation, we shall forget about the indices and just write $Q \in \mathbb{D}_k$ and refer to $Q$ as a cube (or pseudo-cube) of generation $k$. Such cubes satisfy:

\begin{enumerate}[(i)]

\item $\partial \Omega =\bigcup_{j} Q_{j}^k \,\,$ for any $k \in \mathbb Z$.

\item If $m > k$ then either $Q_{i}^{m}\subseteq Q_{j}^{k}$ or
$Q_{i}^{m}\cap  Q_{j}^{k}=\emptyset$.

\item $Q_i^m \cap Q_j^m=\emptyset$ if $i\neq j$. 

\item Each pseudo-cube $Q\in\mathbb{D}_k$ has a ``center'' $x_Q\in \D$ such that
\begin{equation}\label{cube-ball}
\Delta(x_Q,a_02^{-k})\subset Q \subset \Delta(x_Q,2^{-k}).
\end{equation}
\end{enumerate}

If $Q \in \mathbb{D}_{k+\ell}$ for some $\ell\in\N$ and $R$ is the cube of $\mathbb{D}_k$ that contains $Q$ (it is unique by ($ii$), then we say that $R$ is an {\it ancestor} of $Q$, or the {\it parent} of $Q$ if $\ell=1$.

The above definition considers the case where $\partial \Omega$ is unbounded. If the boundary $\partial \Omega$ is bounded, then $\mathbb D_{\partial \Omega} := \bigcup_{k\geq k_0} \mathbb D_k$ where $k_0$ is such that $2^{-k_0-1} < \diam(\Omega) \leq 2^{-k_0}$, it means that there are cubes without ancestors.

If $\mu$ is any doubling measure on $\partial \Omega$ - i.e. if $\mu(2\Delta) \leq C_\mu \mu(\Delta)$ for any boundary ball $\Delta \subset \partial \Omega$ - then we can arrange the following extra property:
\begin{enumerate}[(i)] \addtocounter{enumi}{4}
\item $\mu(\partial Q_i^k) = 0$ for all $i,k$.
\end{enumerate}
Note also that this property makes those sets suitable when weak-$*$ convergences of measures are considered. 

Let us introduce some extra notations. When $E \subset \partial \Omega$ is a set, we define 
$\D_{\partial \Omega}(E):=\set{Q\in\D_{\pom},\, Q\subset E}$.
When $Q\in \mathbb D_{\partial \Omega}$, we write $k(Q)$ for the generation of $Q$ and $\ell(Q)$ for $2^{-k(Q)}$, which is roughly the diameter of $Q$ by \eqref{cube-ball}. The ancestor of $Q$ from $j$ generations ago will be denoted by $Q^{(j)}$. We use $B_Q \subset \R^n$ for the ball $B(x_Q,\ell(Q))$ and $\Delta_Q$ for the boundary ball $\Delta(x_Q,\ell(Q))$ that appears in \eqref{cube-ball}.

\subsubsection*{Whitney decomposition}
We shall also use a covering $\W = \W_\Omega$ of $\Omega$ by dyadic cubes in $\Rn$, and the sets in $\W_\Omega$ are called Whitney cubes. We can construct these Whitney cubes using real dyadic cubes, but we choose to use the dyadic decomposition (of $\Rn$ instead of $\pom$) that we just presented, which is more convenient because we do not need to mention the dimension\footnote{The ``real'' dyadic cubes in $\R^n$ will not satisfy \eqref{cube-ball} with a universal constant $a_0$.}. We define $\W$ as the collection of maximal dyadic cubes $W \in \D_{\R^n}$ such that $20B_W \subset \Omega$. Recall that $B_W$ is the ball centered at $x_W$ with radius $\ell(W)$ in $\Rn$. They have the property that
\[ \Omega = \bigcup_{W\in \W} W\]
and
\[ 20\ell(W) \leq \dist(W,\partial \Omega) < 40\ell(W).\]
We write $W^*$ for $2B_W$, and we observe that $\{W^*\}_{W\in \W}$ covers $\Omega$ and is (uniformly) finitely overlapping.

\subsubsection{Uniform rectifiability}

\begin{definition} \label{defUR}
A $d$-Ahlfors regular (with constant $C_\sigma$) set $\Gamma$ is uniformly $d$-rectifiable if $d$ is an integer, and if there exist constants $\epsilon,M>0$ such that for any $x\in \Gamma$ and $r\in (0,\diam\Gamma)$, there is a $M$-Lipschitz function $f_{x,r}:B(0,r) \subset \R^d \to \R^n$ such that 
\[\sigma(\Gamma \cap B(x,r) \cap f_{x,r}(B(0,r)) \geq \epsilon r^d.\]
We call the triple $(C_\sigma,\epsilon,M)$ the uniform rectifiability constants of $\Gamma$.
\end{definition}

Many equivalent characterizations of uniform rectifiability exist. Let us present the ones that we shall need in this article. We write $\mathcal P_d$ for the collection of affine $d$-plane in $\R^n$. For a $d$-Ahlfors regular set $\Gamma$, we define the bilateral $\beta$-numbers as
\[b\beta_\infty(Q) := \ell(Q)^{-1} \inf_{P\in \mathcal P_d} \Big( \sup_{y \in \Gamma \cap 2B_Q} \dist(y,P) + \sup_{y \in P \cap 2B_Q} \dist(y,\Gamma) \Big).\]

\begin{theorem}[BWGL] \label{ThUR}
Let $\Gamma$ be $d$-Ahlfors regular. Then $\Gamma$ is uniformly rectifiable if and only if, for any $\epsilon>0$, there exists  $C_\epsilon>0$ such that for any $Q_0\in \D_{\Gamma}$,
\begin{equation} \label{BWGL}
\sum_{Q\in \D_{\Gamma}(Q_0) \atop \beta_\infty(Q) >\epsilon} \sigma(Q) \leq C_\epsilon \sigma(Q_0).
\end{equation}
\end{theorem}

\bp
The non-dyadic version, that is, when \eqref{BWGL} is written in terms of Carleson measures -  can be found in \cite[Theorem 2.4]{DS2}. Since any ball  centered at the boundary $B(x,r)$ is included in $2B_Q$ for some $Q\in \D_{\partial \Omega}$ that satisfies $\ell(Q) \leq 4r$, the non-dyadic version of \eqref{BWGL} is equivalent to its dyadic version.
\ep

\subsection{Regularity of solutions}

We start with the definition of the Green function. 
For a uniformly elliptic operator $L=-\diver A\nabla$, one can construct the Green function $G$ on a uniform domain (cf. \cite{KenigB}, and \cite{DFM22} for the higher-codimension case).

\begin{definition}[The Green function]\label{defGreen}
Let $\Omega\subset\R^{n}$ be a uniform domain with $d$-Ahlfors regular boundary. Let $L=-\diver A\nabla$ be a uniformly elliptic operator in $\om$. Then there exists a unique Green function $G_L(X,Y): \Omega \times \Omega  \to \R\cup\set{\infty}$ 
with the following properties: 
$G_L(\cdot,Y)\in W(\Omega\setminus \{Y\})\cap C(\overline{\Omega}\setminus\{Y\})$, 
$G_L(\cdot,Y)\big|_{\pom}\equiv 0$ for any $Y\in\Omega$, 
and $L G_L(\cdot,Y)=\delta_Y$ in the weak sense in $\Omega$, that is,
\begin{equation*}\label{Greendef}
    \iint_\Omega A(X)\,\nabla_X G_{L}(X,Y) \cdot\nabla\vp(X)\, dX=\vp(Y), \qquad\text{for any }\vp \in C_c^\infty(\Omega).
\end{equation*}
In particular, $G_L(\cdot,Y)$ is a positive weak solution to $L G_L(\cdot,Y)=0$ in $\Omega\setminus\{Y\}$. Given $Y\in\om$, we write $G^Y=G_L(\cdot,Y)$ and call it the Green function with pole at $Y$.
\end{definition}

We focus on regularity of solutions to equations with coefficients that satisfy the condition \eqref{DAbdd}. The following result can be found in \cite[Lemma 4.47]{HMT}.

\begin{lemma}\label{lem.intpwest}
    Let $B=B(X_0,r)$ be a ball in $\Rn$. Let $\A$ be an elliptic matrix in $B$ that satisfies 
    \[
\norm{\nabla \A}_{L^\infty(B)}\le \frac{M}{r} \text{  for some }M>0.
\]
 Then there exists some $C_M>0$ depending only on $n$, $M$ and ellipticity such that for every $u\in W^{1,2}(B)\cap L^\infty(B)$, $u\ge0$ that verifies $\diver(\A\nabla u)=0$ in $B$ in the weak sense, we have
 \begin{equation}\label{eq.HMT1}
     \sup_{X\in\frac12 B}\abs{\nabla u(X)}\le \frac{C_M}{r}\inf_{X\in\frac12 B}u(X),
 \end{equation}
 and 
 \begin{equation}\label{eq.HMT2}
     \iint_{\frac14 B}\abs{\nabla^2u(X)}^2dX\le \frac{C_M}{r^2}\iint_{\frac12 B}\abs{\nabla u(X)}^2dX.
 \end{equation}
\end{lemma}

Lemma \ref{lem.intpwest} can be adapted to our setting and immediately yields the following two results. 

\begin{lemma} \label{lemDNu/u<1}
Let $\om\subset\Rn$ be a domain with $d$-Ahlfors regular boundary $\pom$. Take $\beta>0$ and let $L=-\diver\br{D_\beta^{d+1-n}\A\nabla}$ be uniformly elliptic. Assume that $D_\beta |\nabla \A| \leq M$ for some $M>0$. Then for any non-negative weak solution $u$ to $Lu=0$ in $\om$, there holds
\begin{equation}\label{eq.gradupwbd}
    \abs{\nabla u(X)}\le \frac{C u(X)}{\delta_{\partial \Omega}(X)} \quad \text{for all } X\in\om,
\end{equation}
where the constant $C>0$ depends only on $d$, $n$, $M$, $\beta$, and $C_\sigma$, and the elliptic constant of $L$.
\end{lemma}

\begin{lemma} \label{lemuinW22}
Let $\om$ and $L$ be as in Lemma \ref{lemDNu/u<1}. Then for any non-negative weak solution $u$ to $Lu=0$ in $\om$, there holds
\begin{equation}\label{eq.grad2ubd}
    \iint_{B}\abs{\nabla^2u(X)}^2dX\le \frac{C}{r_B^2}\iint_{2B}\abs{\nabla u(X)}^2dX \quad \text{for any } B \text{ that satisfies }6B\subset\om,
\end{equation}
where $r_B$ is the radius of the ball $B$, and the constant $C>0$ depends only on $d$, $n$, $M$, $\beta$, and $C_\sigma$, the elliptic constant of $L$.
In particular, $u\in W^{2,2}_{loc}(\Omega)$.
\end{lemma}
\medskip

\noindent {\em Proof of Lemmas \ref{lemDNu/u<1} and \ref{lemuinW22}: }
Fix any $X_0\in\om$, and let $B_0=B(X_0,\delta_{\partial \Omega}(X_0)/2)$. When $d=n-1$, Lemma \ref{lem.intpwest} immediately gives 
    $\abs{\nabla u(X_0)}\le C u(X_0)/\delta(X_0)$ and \eqref{eq.grad2ubd},
as desired. When $d<n-1$, $u$ is a weak solution of $-\diver\wt A\nabla u=0$ in $B_0$ with $\wt A(X)=\br{\frac{D_\beta(X)}{D_\beta(X_0)}}^{d+1-n}\A$. Observe that $\wt A$ is elliptic (in the usual sense) in $B_0$ and satisfies 
\begin{multline}\label{eq.rescaleA}
    \abs{\nabla \wt A} \leq C_{d,n,\|\A\|_\infty} \frac{D_\beta^{d-n}}{D_\beta^{d+1-n}(X_0)}(\abs{\nabla D_\beta} + D_\beta \abs{\nabla \A}) \\
    \le C_{d,n,\beta,C_\sigma,\|\A\|_\infty,M} \delta_{\partial \Omega}(X_0)^{-1} \quad \text{for }X\in B_0.
\end{multline}
by \eqref{Db=dist}, and since it is easy to check that $|\nabla D_\beta| \leq D_\beta^{-\beta-1}D_{\beta+1}^{\beta+1} \leq C_{\beta,C_\sigma}$. So we can apply Lemma \ref{lem.intpwest} to obtain \eqref{eq.gradupwbd}. 
The proof of \eqref{eq.grad2ubd} is entirely similar.
\ep

Moreover, one can show that under the assumptions of Lemma \ref{lemDNu/u<1}, $\nabla u$ is H\"older continuous. We formulate this result in a way that is convenient for our application later. 

\begin{lemma}\label{lemDNu/uHolder}
    Let $\om$ and $L$ be as in Lemma \ref{lemDNu/u<1}. Let $W\in\W_{\om}$ and $u$ be a positive solution to $Lu=0$ in $W^*$. Then $u\in C^{1,\gamma}(\overline{W})$ for some $\gamma\in(0,1)$. Moreover, for any $X$, $Y\in W$ with $\abs{X-Y}<\ell(W)/2$, we have 
    \begin{equation}\label{eq.DNu/uHolder}
    \abs{\frac{D_\beta(X)\abs{\nabla u(X)}}{u(X)}-\frac{D_\beta(Y)\abs{\nabla u(Y)}}{u(Y)}}\le C\br{\frac{\abs{X-Y}}{\ell(W)}}^\gamma,    
    \end{equation}
    where $C>0$ depends only on dimension, $M$, $\beta$, and the elliptic constant of $L$.
\end{lemma}

\bp
We first observe that we only need to consider the case when $d=n-1$ because the case when $0<d<n-1$ can be treated similarly (with modifications that are the same as in the proof of Lemma \ref{lemDNu/u<1} and \ref{lemuinW22}). 

We claim that for any $X$, $Y\in W$ with $\abs{X-Y}<\ell(W)/2$,
\begin{equation}\label{eq.NuHolder}
    \abs{\nabla u(X)-\nabla u(Y)}\lesssim \norm{u}_{L^\infty(W^*)}\frac{\abs{X-Y}^\gamma}{\ell(W)^{1+\gamma}}
\end{equation}
for some $\gamma\in(0,1)$. 
This result is probably proven, but the authors could not pinpoint a precise reference. We know that 
in \cite{DK} Theorem 1.5, it was shown that under the weaker assumption on the operator ($L$ have Dini mean oscillation), $u\in C^1(\overline{W})$, and an upper bound of the modulus of continuity of $\nabla u$ is given. Following the proof of  \cite{DK} Theorem 1.5, we can get that for $X$, $Y\in W$ with $r=\abs{X-Y}<\ell(W)/2$,
\[
\abs{\nabla u(X)-\nabla u(Y)}\lesssim\br{\frac{r}{\ell(W)}}^\gamma\fint_{W^*}\abs{\nabla u}dZ+\frac{r}{\ell(W)}\norm{\nabla u}_{L^\infty(W^*)}
\]
for some $\gamma\in(0,1)$, which yields \eqref{eq.NuHolder} by using \eqref{eq.gradupwbd}.

The estimate \eqref{eq.DNu/uHolder} follows immediately from \eqref{eq.NuHolder} as we can write 
\begin{multline*}
    \abs{\frac{D_\beta(X)\abs{\nabla u(X)}}{u(X)}-\frac{D_\beta(Y)\abs{\nabla u(Y)}}{u(Y)}}\\
    \le \frac{D_\beta(X)}{u(X)}\abs{\nabla u(X)-\nabla u(Y)}+\abs{\frac{D_\beta(X)}{u(X)}-\frac{D_\beta(Y)}{u(Y)}}\abs{\nabla u(Y)}.
\end{multline*}
The inequality \eqref{eq.NuHolder} gives the desired upper bound for the first term on the right-hand side. For the second term, since $\abs{\nabla D_\beta}\le C$, we have 
\[
    \abs{\frac{D_\beta(X)}{u(X)}-\frac{D_\beta(Y)}{u(Y)}}\abs{\nabla u(Y)}
    \le 
    \frac{D_\beta(X)\norm{\nabla u}_{L^\infty(W)}+C\norm{u}_{L^\infty(W)}}{\norm{u}^2_{L^\infty(W)}}\abs{\nabla u(Y)}\abs{X-Y},
\]
which is bounded above by $C\frac{\abs{X-Y}}{\ell(W)}$ thanks to \eqref{eq.gradupwbd}.
\ep

\section{Estimates on $D_\beta$ and its derivatives}\label{SecD}

We want to compare $D_\beta$ and its derivatives with a quantity which is easier to deal with. More precisely, we want to replace the integral on $\partial \Omega$ in \eqref{defDbeta} with the integral on a plane, and say that the error is controlled. Note that similar  computations have been done in various contexts (see for instance \cite{DFMdahlberg}, \cite{DEM}, \cite{FenAinfty}, \cite{DLM1}).

We shall follow the presentation of \cite{FenAinfty}. For this, let us introduce some notations first. Let $\partial \Omega$ be $d$-Ahlfors regular and $\sigma$ be as in \eqref{defADR}. For $\beta>0$, let $R_\beta(X)$ be the integral defined on  $\R^n \setminus \partial \Omega$ as
\[R_\beta(X) := \int_{\partial \Omega} |X-y|^{-d-\beta} \, d\sigma(y).\]
We then set 
\[c_\beta := \int_{\R^d} (1+|y|^2)^{-\frac{d+\beta}{2}} \, dy.\]
Given a multiindex $\kappa=(\kappa_1,\dots,\kappa_n)\in\N^n$, we define $\abs{\kappa}:=\sum_{i=1}^n\kappa_i$, and 
\[
\partial^\kappa :=\frac{\partial^{\abs{\kappa}}}{\partial_1^{\kappa_1}\dots\partial_n^{\kappa_n}}.
\] 

\begin{theorem} \label{ThDb}
Let $\partial \Omega$ be uniformly $d$-rectifiable set in $\Rn$. Then there exists $C_\kappa>0$  such that 
\begin{equation} \label{DkRb<}
|\partial^\kappa R_\beta| \leq C_\kappa \delta_{\partial \Omega}^{-\beta - \abs{\kappa}}.
\end{equation}
In addition, there exists $C>0$ and a function $a_\beta \in CM_{\Omega}(M)$ for some $M>0$, such that for any $X\in \R^n \setminus \Omega$, there exists a constant $c_X>0$ and an affine $d$-plane $P_X$ that satisfy
\begin{equation} \label{cX=1}
C^{-1}  \leq c_X \leq C, 
\end{equation}
\begin{equation} \label{DbX=dist}
C^{-1} \delta_{\partial \Omega}(X) \leq D_{\beta,X}(X) \leq C \delta_{\partial \Omega}(X), 
\end{equation}
where $D_{\beta,X}(Y):= (c_\beta)^{-1/\beta} \dist(Y,P_X)$, and
\begin{equation} \label{Rb=RbX}
\left| \partial^\kappa R_\beta(X) - c_X \partial^\kappa R_{\beta,X}(X) \right| \leq C_\kappa \delta_{\partial \Omega}(X)^{-\beta-\abs{\kappa}} a_\beta(X),
\end{equation}
where $R_{\beta,X}(Y) := c_\beta \dist(Y,P_X)^{-\beta}$.  One possible choice of $c_X$ is  
\begin{equation} \label{cX=sX}
c_X := \frac{c_1}{c_{1/2}^{2}} \frac{D_1(X)}{D_{1/2}(X)}.
\end{equation}
The constant $C$ depends on $\beta$ and the Ahlfors-regular constant in \eqref{defADR}, $C_\kappa$ depends also on $\kappa$, $M$ depends on $\beta$ and the uniformly rectifiable constants of $\partial \Omega$.
\end{theorem}

\begin{remark}
\begin{itemize}
\item $P_X$ is morally the plane that best approximates $\partial \Omega$ in $B(X,100\delta_{\partial \Omega}(X))$.
\item The function $a_\beta$ is morally a geometric sum of Tolsa's $\alpha$-numbers, which are numbers that can be used to characterize the uniform rectifiability (see \cite[Theorem 1.2]{Tol09}).
\end{itemize}
\end{remark}

\bp
The uniform bound \eqref{DkRb<} is fairly simply, since $|\partial^\kappa R_\beta| \leq C_{d,\beta,k} R_{\beta+\abs{\kappa}}$ and the bound $R_\beta \leq C\delta_{\partial \Omega}^{-\beta}$ is implied by \eqref{Db=dist}.

The equivalence \eqref{DbX=dist} is a consequence of the choice of $P_X$ in \cite{FenAinfty}, and is stated as (3.9)--(3.10) in \cite{FenAinfty}\footnote{The plane $P_X$ depends only on the Whitney cube that contains $X$, hence the notation $P_Q$ used in \cite{FenAinfty}}.
For \eqref{Rb=RbX}, the result can be found in \cite{FenAinfty} Lemma 3.31 and (3.30) when $\partial^\kappa$ is a derivative of order 0 or 1, but the proof easily extends to derivatives of higher order. Theorem 2.1 in \cite{DEM} provides another similar case involving second derivatives.

The fact that $a_\beta \in CM$ when $\partial \Omega$ is uniformly rectifiable is (3.51) in \cite{FenAinfty} or similarly \cite[Lemma 5.89]{DFMdahlberg}. 

Being able to choose $c_X$ as in \eqref{cX=sX} is a consequence of (3.63) in \cite{FenAinfty}.
\ep

\begin{corollary} \label{CorDb}
Under the same conditions and notation as in Theorem \ref{ThDb}, we have
\begin{equation} \label{DkDb<}
|\partial^\kappa D_\beta| \leq C_\kappa \delta_{\partial \Omega}^{1 - \abs{\kappa}},
\end{equation}
for any $\nu\in \R$,
\begin{equation} \label{Dnu=DnuX}
|(D_\beta)^\nu(X) - (c_X)^{-\nu/\beta} D_{\beta,X}^\nu(X)| \leq C_\nu \delta_{\partial \Omega}(X)^{\nu} a_\beta(X),
\end{equation}
and
\begin{equation} \label{Db=DbX}
\left| \partial^\kappa D_\beta(X) - (c_X)^{-1/\beta} \partial^\kappa D_{\beta,X}(X) \right| \leq C_\kappa \delta_{\partial \Omega}(X)^{1-\abs{\kappa}} a_\beta(X).
\end{equation}
\end{corollary}

\bp
Recall that $D_\beta = R_\beta^{-1/\beta}$, so we write $\partial^\kappa D_\beta$ as a sum of products of powers of $R_\beta$ and derivatives of $R_\beta$. The bound \eqref{DkDb<} is then a consequence of \eqref{Db=dist} and \eqref{DkRb<}.

The estimate \eqref{Dnu=DnuX} is equivalent to 
\begin{equation} \label{Dnu=DnuX2}
|(R_\beta)^{-\nu/\beta}(X) - (c_X)^{-\nu/\beta} R_{\beta,X}^{-\nu/\beta}(X)| \leq C_\nu \delta_{\partial \Omega}(X)^{\nu} a_\beta(X).
\end{equation}
To prove the above bound, we apply the Mean Value Theorem to the function $z\mapsto z^{-\nu/\beta}$ and we get
\begin{multline*}
|(R_\beta)^{-\nu/\beta}(X) - (c_X)^{-\nu/\beta} R_{\beta,X}^{-\nu/\beta}(X)| \\
\leq \frac{|\nu|}{\beta} \max\{(R_\beta)^{-\nu/\beta-1}(X), c_X^{-\nu/\beta-1}R_{\beta,X}^{-\nu/\beta-1}(X)\} |R_\beta(X) - c_X R_{\beta,X}(X)|  \\
\leq C_{\beta,\nu} \delta_{\partial \Omega}(X)^{\nu+\beta} |R_\beta(X) - c_X R_{\beta,X}(X)| 
\end{multline*}
thanks to \eqref{Db=dist}, \eqref{cX=1}, and \eqref{DbX=dist}. We conclude \eqref{Dnu=DnuX2} by invoking \eqref{Rb=RbX}.

Recall that $\partial^\kappa D_\beta$ is a sum of products of powers of $D_\beta$ and derivatives of $R_\beta$. So, \eqref{Db=DbX} is a straightforward consequence of \eqref{Dnu=DnuX}, \eqref{Rb=RbX}, and the two following stability results (whose proofs are immediate):
\begin{enumerate}[(i)]
\item {\em Stability under sum.} If $f_1$, $f_2$, $g_1$ and $g_2$ are 4 functions on $\Omega$ such that, for $i\in \{1,2\}$,  
\[|f_i(X) - g_i(X)| \leq C\delta_{\partial \Omega}(X)^\nu a_\beta(X),\]
then
\[|(f_1+f_2)(X) - (g_1+g_2)(X)| \leq 2C\delta_{\partial \Omega}(X)^\nu a_\beta(X).\]
\item {\em Stability under product.} If  If $f_1$, $f_2$, $g_1$, $g_2$ are 4 functions on $\Omega$ and $\nu_1$, $\nu_2$ are two real numbers such that, for $i\in \{1,2\}$,
\[|f_i(X)| + |g_i(X)| \leq C\delta_{\partial \Omega}(X)^{\nu_i}\]
and
\[|f_i(X) - g_i(X)| \leq C\delta_{\partial \Omega}(X)^{\nu_i} a_\beta(X),\]
then
 and
\[|(f_1f_2)(X) - (g_1g_2)(X)| \leq 2C^2\delta_{\partial \Omega}(X)^{\nu_1+\nu_2} a_\beta(X).\]
\end{enumerate}
The corollary follows. \ep

\begin{corollary} \label{CorDb2}
Under the same conditions and notation as as Theorem \ref{ThDb}, for any $X\in\Rn\setminus\pom$, we have 
\begin{equation} \label{NNDb2isCM}
\big|\nabla\br{ |\nabla D_\beta(X)|^2}\big| \leq C \delta_{\partial \Omega}(X)^{-1} a_\beta(X),
\end{equation}
and 
\begin{equation} \label{LbDbisCM}
|\div[D_\beta^{d+1-n}(X) \nabla D_\beta(X)]| \leq C \delta_{\partial \Omega}(X)^{d-n} a_\beta(X),
\end{equation}
and therefore,
\[ D_\beta \big|\nabla \br{|\nabla D_\beta|^2}\big| + D_\beta^{n-d} |\div(D_\beta^{d+1-n} \nabla D_\beta)| \in CM_\Omega(CM).\]

If $d=n-1$, then we even have $D_\beta \abs{\nabla^2 D_\beta} \in CM_\om(CM)$, or more precisely,
\begin{equation} \label{N2DbisCM}
\big|\nabla^2 D_\beta(X)\big| \leq C \delta_{\partial \Omega}(X)^{-1} a_\beta(X).
\end{equation}
\end{corollary}

\bp When $d=n-1$, $D_{\beta,X}$ is - up to a constant - the distance to a hyperplane, so $\nabla^2 D_{\beta,X} = 0$. The bound \eqref{N2DbisCM} follows from \eqref{Db=DbX} (with $\abs{\kappa}=2$). 

\medskip

For the other bounds, observe that 
\[\nabla \br{|\nabla D_\beta|^2} = 2 \sum_{k=1}^n (\nabla \partial_k D_\beta) (\partial_k D_\beta),\]
and
\[ \div (D_\beta^{d+1-n} \nabla D_\beta)  = (d+1-n) D_\beta^{d-n} \sum_{k=1}^n (\partial_k D_\beta)^2 + D^{d+1-n} \Delta D_\beta.\]
 that is, both quantities are sum of products of derivatives of $D_\beta$ (and powers of $D_\beta$). So \eqref{Dnu=DnuX}, \eqref{Db=DbX}, and the stability under sum and product mentioned in the proof of Corollary \ref{CorDb} give that 
 \[\abs{\nabla \br{|\nabla D_\beta(X)|^2} - \nabla \br{|\nabla D_{\beta,X}(X)|^2}} \leq C \delta_{\partial \Omega}(X)^{-1} a_\beta(X)\]
and
\[|\div[D_\beta^{d+1-n}(X) \nabla D_\beta(X)] - \div[D_{\beta,X}^{d+1-n}(X) \nabla D_{\beta,X}(X)]| \leq C \delta_{\partial \Omega}(X)^{d-n} a_\beta(X).\]
The corollary will be proven if we can show that 
\[\nabla (|\nabla D_{\beta,X}(X)|^2)= 0 \quad \text{ and } \quad  \div[D_{\beta,X}^{d+1-n}(X) \nabla D_{\beta,X}(X)] = 0.\] 
For the first equality, recall that $D_{\beta,X}$ is - up to the multiplicative constant $c:=c_\beta^{-1/\beta}$ - a distance to a $d$-plane, so for $Y$ that does not touch $P_X$, we have $|\nabla D_{\beta,X}(Y)| = c$ and $\nabla (|\nabla D_{\beta,X}(Y)|^2) = 0$, which proves $\nabla |\nabla (D_{\beta,X}(X)|^2) = 0$. 
For the second equality, we can assume that $P_X = \R^d \times \{0\}$ and write a point in $\R^n$ as $Y=(y,t) \in \R^d \times \R^{n-d}$. Thus $D_{\beta,X}(Y) = c |t|$ and $\div[D_{\beta,X}^{d+1-n}(X) \nabla D_{\beta,X}(X)] = c^{d+2-n} \div[ |t|^{d+1-n} \nabla |t|] = 0$ as desired. The corollary follows.
\ep

\begin{corollary} \label{CorDb3}
Under the same conditions and notation as as Theorem \ref{ThDb}, there exists $c = c(C_\sigma,\beta)>0$ such that for any $X\in\Rn\setminus\pom$,
\begin{equation} \label{1NDb<1isCM}
\1_{\{|\nabla D_\beta|< c\}}(X) \leq C a_\beta(X)
\end{equation}
and
\begin{equation} \label{NNDbisCM}
|\nabla |\nabla D_\beta(X)|| \leq C \delta_{\partial \Omega}(X)^{-1} a_\beta(X).
\end{equation}
Therefore,
\[ \1_{\{|\nabla D_\beta| < c\}} + D_\beta \big|\nabla |\nabla D_\beta|\big| \in CM_\Omega(CM).\]
\end{corollary}

\bp
Assume first that \eqref{1NDb<1isCM} is true and let us prove \eqref{NNDbisCM}. Thanks to \eqref{DkDb<}, we have that $|\nabla|\nabla D_\beta|| \leq C\delta_{\partial \Omega}^{-1}$. Together with the fact that 
\[ \big|\nabla |\nabla D_\beta|\big| = \dfrac{\big|\nabla |\nabla D_\beta|^2\big|}{|\nabla D_\beta| },\]
we have
\[ \big|\nabla |\nabla D_\beta(X)|\big| \leq C\1_{\{|\nabla D_\beta(X)| < c\}} \delta_{\partial \Omega}^{-1}(X) + c^{-1} \big|\nabla |\nabla D_\beta(X)|^2\big| \leq C'\delta_{\partial \Omega}(X)^{-1} a_\beta(X)\]
by \eqref{NNDb2isCM} and \eqref{1NDb<1isCM}.

It remains to show \eqref{1NDb<1isCM}. Recall that $c_X \approx 1$ by \eqref{cX=1}, and that 
$D_{\beta,X}(Y)=c_\beta^{-1/\beta}\dist(Y,P_X)$. So 
$|\nabla D_{\beta,X}(X)| = (c_\beta)^{-1/\beta}$, and thus, there exists $c>0$ such that, for any $X\in \Rn\setminus\pom$, $|(c_X)^{-1/\beta}\nabla D_{\beta,X}(X)| \geq 2c$. Then we have
\begin{multline*}
\1_{\{|\nabla D_\beta| < c\}}(X) \leq \frac1c \big| |\nabla D_\beta(X)| - |(c_X)^{-1/\beta}\nabla D_{\beta,X}(X)| \big| \\
\leq \frac1c \big| \nabla D_\beta(X) - (c_X)^{-1/\beta}\nabla D_{\beta,X}(X) \big| 
\leq C a_\beta(X)
\end{multline*}
by the triangle inequality and then \eqref{Db=DbX}. The corollary follows.
\ep

\section{The Caccioppoli-type estimate}\label{SecCaccio}

The aim of the section is to deduce Carleson estimates on the second derivatives of $\ln(u/D)$ from the first ones. 

\begin{theorem} \label{ThCaccio}
Let $d\in [0,n)$, $\Omega \subset \R^n$ be a domain with $d$-Ahlfors regular boundaries, and $L:=-\diver(D^{d+1-n} \A \nabla)$ be uniformly elliptic with respect to a weight $D^{d+1-n}$ that satisfies $D\in W_{\loc}^{2,\infty}(\om)$, 
\begin{equation} \label{ThCaccioK}
K^{-1} \dist(X,\partial \Omega) \leq D(X) \leq K\dist(X,\partial \Omega), \  \text{ and } \ |\nabla D(X)| \leq K \qquad \text{ for } X \in \Omega.
\end{equation}
Let $W\in \W_\Omega$, and assume that $D\abs{\nabla \A}$ is bounded on $W^*$. Then for any non-negative non-trivial weak solution $u$ to $Lu=0$ in $W^*$, we have
\begin{multline}\label{cacci.codim1}
    \iint_{W} \left|\nabla^2 \ln \Big( \frac{u}{D} \Big) \right|^2 D^{d+4-n} dX \leq C \Big\{\iint_{W^*} \left|\nabla \ln \Big( \frac{u}{D} \Big) \right|^2 D^{d+2-n} \, dX \\
    + \iint_{W^*} |D^{n-d}\diver(D^{d+1-n} \nabla D)|^2  \, D^{d+2-n} \, dX + \iint_{W^*} |\nabla \A|^2 D^{d+2-n} \, dX \Big\},
\end{multline}
where $C$ depends only on $d$, $n$, the ellipticity constant of $L$, $\|D\nabla A\|_{L^\infty(W^*)}$, and the constant $K$ in \eqref{ThCaccioK}.
\end{theorem}

\begin{corollary} \label{CorCaccio}
Let $\Omega \subset \R^n$ be such that $\partial \Omega$ is uniformly $d$-rectifiable. Take $\beta >0$ and define $D_\beta$ as in \eqref{defDbeta}. Let $L=-\diver\br{D_\beta^{d+1-n}\A\nabla}$ be uniformly elliptic. Assume either
\begin{enumerate}
    \item $d=n-1$, $\Omega$ is uniform, and $\A$ satisfies the DKP condition, or
    \item $0<d<n-1$ and $\A=I$.
\end{enumerate}
Then 
for any $x\in \partial \Omega$, any $r\in (0,\diam \partial \Omega)$ and any positive weak solution $u$ to $Lu=0$ in $B(x,2r)$ such that $u=0$ on $\partial \Omega \cap B(x,2r)$, we have
\[\iint_{B(x,r)} \left|\nabla^2 \ln \Big( \frac{u}{D_\beta} \Big) \right|^2 D_\beta^{d+4-n} dX \leq C r^d,\]
where $C>0$ depends on $d$, $n$, $\beta$, the elliptic, boundedness and DKP constant of $L$, and uniform rectifiability constants of $\partial \Omega$. 
\end{corollary}

\noindent {\em Proof of Corollary \ref{CorCaccio}.} The fact that $D_\beta$ verifies \eqref{ThCaccioK} was already stated in \eqref{Db=dist} and \eqref{DkDb<}, and the fact that $\|D_\beta \nabla \A\|_\infty < \infty$ comes from the definition of DKP-operators. Let $\C(x,r)$ be the subcollection of elements in $\W$ that intersect $B(x,r)$. Note that by construction, $W^* \subset B(x,3r/2)$ when $W \in \C(x,r)$. Theorem \ref{ThCaccio} entails
\begin{multline*}
\iint_{B(x,r)} \left|\nabla^2 \ln \Big( \frac{u}{D_\beta} \Big) \right|^2 D_\beta^{d+4-n} dX \leq \sum_{W\in \C(x,r)} \iint_{W} \left|\nabla^2 \ln \Big( \frac{u}{D_\beta} \Big) \right|^2 D_\beta^{d+4-n} dX \\
\lesssim  \sum_{W \in \C(x,r)}  \iint_{W^*} \Big\{ \left|D_\beta \nabla \ln \Big( \frac{u}{D_\beta} \Big) \right|^2 + |D_\beta^{n-d}\diver(D_\beta^{d+1-n} \nabla D_\beta)|^2 + |D_\beta \nabla \A|^2 \Big\}  D_\beta^{d-n} dX \\
\lesssim \iint_{B(x,3r/2)} \Big\{ \left|D_\beta \nabla \ln \Big( \frac{u}{D_\beta} \Big) \right|^2 + |D_\beta^{n-d}\diver(D_\beta^{d+1-n} \nabla D_\beta)|^2 + |D_\beta \nabla \A|^2 \Big\} D_\beta^{d-n} dX \\
=: T_1 + T_2 + T_3,
\end{multline*}
since the $W^*$ are finitely overlapping. The last term $T_3$ is bounded by $\sigma(Q)$ by definition of $L$ being a DKP-operator. The bounds on $T_1$ and $T_2$ come from Theorem \ref{ThNlnu/D} and Corollary \ref{CorDb2} respectively. 
\ep

\medskip

\noindent {\em Proof of Theorem \ref{ThCaccio}.} We first recall the following facts: for any $X\in \Omega$, we have 
\begin{equation} \label{ND<1}
|\nabla D(X)| \leq K,
\end{equation} 
\begin{equation} \label{DNu/u<1}
\frac{|D\nabla u(X)|}{u(X)}  \leq C
\end{equation} 
and
\begin{equation} \label{DNlnu/D<1}
\left|D\nabla \ln\Big(\frac{u}{D}\Big)\right| = D\left| \frac{\nabla u}{u} - \frac{\nabla D}{D}\right| \leq C.
\end{equation} 
The bound \eqref{ND<1} holds by assumption, the bound \eqref{DNu/u<1} comes from Lemma \ref{lemDNu/u<1}, and the last bound  \eqref{DNlnu/D<1} is just a combination of the two previous ones.

\medskip
Construct a cut-off function $\phi \in C_0^\infty(W^*)$ such that $0 \leq \phi \leq 1$, $\phi \equiv 1$ on $W$, and $|\nabla \phi| \leq 100/ \dist(.,\partial \Omega) \lesssim 1/D.$ It is easy to check that the last condition together with \eqref{ND<1} implies
\begin{equation} \label{NDjphi2<1}
|\nabla (D^j \phi^2)| \lesssim D^{j-1} \phi
\end{equation} 
for any $j\in \mathbb N$ (with a constant that depends on $j$).

Let $\partial_k$ be the $k^{th}$ coordinate of the gradient. Let us introduce $T$, $U_1$, $U_2$, $U_3$ as the integral terms:
\[\begin{split}
T & := \iint_{\Omega} \A \nabla \partial _k \ln \Big( \frac{u}{D} \Big) \cdot \nabla \partial _k \ln \Big( \frac{u}{D} \Big) \,  D^{4+d-n} \phi^2\,  dX \\
U_1 & := \iint_{W^*} \left|\nabla \ln \Big( \frac{u}{D} \Big) \right|^2 D^{2+d-n} \, dX \\
U_2 & := \iint_{W^*} |\diver(D^{d+1-n} \nabla D)|^2  \, D^{n-d} \, dX \\
U_3 & := \iint_{W^*} |\nabla \A|^2 D^{2+d-n} \, dX.
\end{split}\]
Note that $\partial _k \ln \big( \frac{u}{D} \big)\in W_{\loc}^{1,2}(W^*)$ due to Lemma \ref{lemuinW22} and the assumption that $D\in W_{\loc}^{2,\infty}(\om)$. So $T$ is finite. Therefore, thanks to the uniform ellipticity and boundedness of $\A$ given in \eqref{defelliptic}--\eqref{defbounded}, the bound \eqref{cacci.codim1} will be proved once we establish
\begin{equation} \label{boundT<U}
T \lesssim T^{1/2} (U_1 + U_2 + U_3)^{1/2}.
\end{equation}
For the sake of precision, we need to introduce a smooth approximation of $\partial _k \ln \big( \frac{u}{D} \big)$. We choose $v_k$ to be a smooth function that satisfies 
\begin{equation} \label{Tandvk}
T \leq 2 \iint_\Omega \A \nabla \partial _k \ln \Big( \frac{u}{D} \Big) \cdot \nabla v_k \, D^{4+d-n} \phi^2 \, dX \ \text{ and } \ \iint_{\Omega} |\nabla v_k|^2 \, D^{4+d-n} \phi^2 \, dX \leq C T,
\end{equation}
where $C$ depends only on the constant in \eqref{defbounded}.

\medskip

We are ready for the main argument. We use an integration by part to move the $\partial_k$ away from $\ln \big( \frac{u}{D} \big)$, and we get
\begin{multline*}
\frac12 T \leq   \iint_\Omega \A \nabla \partial _k \ln \Big( \frac{u}{D} \Big) \cdot \nabla v_k \, D^{4+d-n} \phi^2 \, dX \\
 =  - \iint_\Omega \partial_k (\A) \nabla \ln \Big( \frac{u}{D} \Big) \cdot \nabla v_k \, D^{4+d-n} \phi^2 \, dX 
 - \iint_\Omega \A \nabla \ln \Big( \frac{u}{D} \Big) \cdot \nabla v_k \, \partial_k[D^{4+d-n} \phi^2] \, dX
 \\
 - \iint_\Omega \A \nabla \ln \Big( \frac{u}{D} \Big) \cdot \nabla \partial_k v_k \, D^{4+d-n} \phi^2 \, dX =: T_1 + T_2 + T_3
\end{multline*}
We develop $T_3$ further: we integrate by part again to move the gradient again from $\partial_k v_k$ and we get
\begin{multline*}
T_3  =  \iint_\Omega \A \nabla \ln \Big( \frac{u}{D} \Big) \cdot \nabla[ D^{2} \phi^2]  \, D^{2+d-n} \partial_k v_k \,  \, dX 
+ \iint_\Omega \div\Big[D^{2+d-n} \A \nabla \ln \Big( \frac{u}{D} \Big)\Big] \,  \partial_k v_k \, D^2 \phi^2 \, dX \\
 =: T_4 + T_5.
\end{multline*}
The terms $T_2$ and $T_4$ are treated similarly: we use \eqref{NDjphi2<1}, the Cauchy-Schwarz inequality, the boundedness of $\A$, and \eqref{Tandvk} to obtain
\begin{multline*}
|T_2| + |T_4| \lesssim \left(\iint_\Omega |\nabla v_k|^2 D^{4+d-n} \phi^2 \, dX \right)^\frac12 \left( \iint_{W^*} \left|\nabla \ln \Big( \frac{u}{D} \Big)\right|^2 D^{2+d-n} \, dX \right)^\frac12
\lesssim T^{1/2} (U_1)^{1/2}.
\end{multline*}
The term $T_1$ is also easy: by the Cauchy-Schwarz inequality, we have
\[|T_1| \lesssim T^\frac12 \left( \iint_{\Omega} |\nabla \A|^2 \left|D \nabla \ln \Big( \frac{u}{D} \Big)\right|^2 D_\beta^{d+2-n} \phi^2 \, dX \right)^\frac12
\lesssim T^{1/2} (U_3)^{1/2}\]
by \eqref{DNlnu/D<1}. It remains to deal with $T_5$. Note first that this term makes sense since $u\in W^{2,2}_{loc}(W^*)$ by Lemma \ref{lemuinW22} and $\nabla\A \in L^\infty_{loc}(W^*)$ by assumption. Moreover, the function $-\div(D^{d+1-n} \A \nabla u)$ is also well defined in $L^2(W^*)$ and is zero almost everywhere. Let us first remove the $v_k$ by applying the Cauchy-Schwarz inequality
\[|T_5| \lesssim T^\frac12 \left( \iint_{\Omega} \left| \div\Big[D^{2+d-n} \A \nabla \ln \Big( \frac{u}{D} \Big)\Big]\right|^2 D^{n-d} \phi^2 \, dX \right)^\frac12 =: T^{1/2} (T_6)^{1/2}.\]
To continue, we need to separate $u$ from $D$, so we use the fact that $\nabla \ln(u/D) = \frac{\nabla u}{u} - \frac{\nabla D}{D}$, and we obtain
\[T_6 \lesssim \iint_{\Omega} \left| \div(D^{d+1-n} \A \nabla D) \right|^2 D^{n-d} \phi^2 \, dX +  \iint_{\Omega} \left| \div\Big( \frac{D}{u} D^{d+1-n} \A \nabla u \Big) \right|^2 D^{n-d} \phi^2 \, dX =: T_7 + T_8.\]
For the term $T_7$, we still need to remove $\A$, so we write
\[T_7 \lesssim \iint_{\Omega} \left| \div(D^{d+1-n} \nabla D) \right|^2 D^{n-d} \phi^2 \, dX  + \iint_{\Omega} | \nabla \A|^2 |\nabla D|^2 D^{2+d-n} \phi^2 \, dX \lesssim U_2 + U_3\]
by \eqref{ND<1}. As for $T_8$, since $u \in W^{2,2}_{loc}(W^*)$ is a solution to $L u = -\div(D^{d+1-n} \mathcal A\nabla u) = 0$,  we have
\[D^{n-d-1} \left|\div\Big(\frac Du D^{d+1-n} \mathcal A\nabla u\Big)\right| = 
\left|\nabla\Big(\frac{D}{u}\Big) \cdot \mathcal A \nabla u\right|
\lesssim \frac{|D\nabla u|}{|u|}\abs{\nabla\ln\br{\frac{u}{D}}},\]
and hence
\[ T_8 \lesssim \iint_{\Omega} \left|\nabla\ln \Big(\frac{u}{D}\Big)\right|^2 \frac{|D\nabla u|^2}{|u|^2} D^{2+d-n} \phi^2 \, dX \lesssim U_1.
\]
by \eqref{DNu/u<1}. The theorem follows.
\ep

\section{UR implies second derivative estimates for $G$}\label{SecURtoG}

Let $\om$ be a domain in $\Rn$ with uniformly $d$-rectifiable boundary. In this section, we prove that the (appropriately normalized) second derivatives of a positive solution $u$ can be controlled by functions that satisfy the Carleson measure condition in $\om$ when $\pom$ is uniformly rectifiable. 

We start with writing \begin{equation}\label{eq.Dudecomp}
   \nabla u = D\nabla\br{\frac{u}{D}}+\frac{u}{D}\nabla D =u\nabla\ln\br{\frac{u}{D}}+\frac{u}{D}\nabla D,
\end{equation}
where $D$ will be the regularized distance $D_\beta$ in application. 
We collect our computations on the second derivatives of $u$ in the following lemma. 

\begin{lemma}\label{lemD2uptwise}
     Let $u\in W_{\loc}^{2,2}(\om)$ and $D\in W_{\loc}^{2,\infty}(\om)$ be two positive functions that satisfy 
    \begin{equation}\label{cond.uD}
        \frac{D\abs{\nabla u}}{u}\le C\quad \text{and } \abs{\nabla^j D}\le C D^{1-j} \quad\text{for }j=1,2\quad \text{in }\om.
    \end{equation}
    Then the following estimates hold almost everywhere in $\om$.
    \begin{equation}\label{eq.D2u}
        \frac{D^2\abs{\nabla^2u}}{u}\lesssim D\abs{\nabla\ln\br{\frac{u}{D}}}+D^2\abs{\nabla^2\ln\br{\frac{u}{D}}}+D\abs{\nabla^2D},
    \end{equation}
     \begin{equation}\label{eq.D|Du|2}
         \frac{D^3\abs{\nabla\br{\abs{\nabla u}^2}}}{u^2}\lesssim D^2\abs{\nabla^2\ln\br{\frac{u}{D}}}+D\abs{\nabla\ln\br{\frac{u}{D}}}+D\abs{\nabla(\abs{\nabla D}^2)},
    \end{equation}
       where the implicit constants depend on $C$ in \eqref{cond.uD}.
\end{lemma}
While we state the lemma in such general terms, in application, we always take $D=D_\beta$ to be the smooth distance and $u$ to be a positive solution to $Lu=0$ in $\om$, with $L$ satisfies the conditions of Theorem \ref{mainthm}. Therefore, the condition \eqref{cond.uD} is satisfied.
\medskip

\bp
We let the reader check that \eqref{eq.D2u} follows easily from taking derivative of \eqref{eq.Dudecomp}, which gives 
\begin{multline*}
\frac{D^2 |\nabla^2 u|}{u} \leq D \frac{D|\nabla u|}{u} \abs{\nabla \ln\br{\frac uD}} + D^2 \abs{\nabla^2 \ln\br{\frac uD}} + D  \abs{\frac{D}{u} \nabla\br{\frac uD}} |\nabla D| +  D |\nabla^2 D|,  
\end{multline*}
then by observing that $\frac{D}{u} \nabla\br{\frac uD} = \nabla \ln\br{\frac uD}$, and finally from invoking the assumption \eqref{cond.uD}.

To get \eqref{eq.D|Du|2}, let us first recall that $D \abs{\nabla \ln\br{\frac{u}{D}}}$ is bounded. Indeed,
\begin{equation}\label{eq.DlnuD}
    D\abs{\nabla\ln\br{\frac{u}{D}}}=\abs{\frac{D\nabla u}{u}-\nabla D} \leq C,
\end{equation}
because of our assumption \eqref{cond.uD}. Then we differentiate the norm of \eqref{eq.Dudecomp} and we obtain
\begin{multline*}
\frac{D^3 \nabla (|\nabla u|^2)}{u^2} = \frac{D^3}{u^2} \nabla\br{\abs{ u \nabla \ln\br{\frac{u}{D}}}^2} + 2 \frac{D^3}{u^2} \nabla \br{ \frac{u^2}{D} \nabla \ln\br{\frac{u}{D}} \cdot \nabla D} \\
+ \frac{D^3}{u^2} \nabla\br{ \frac{u^2}{D^2} |\nabla D|^2} =: T_1 + T_2 + T_3.
\end{multline*}
The term $T_1$ is bounded as follows
\begin{multline*}
|T_1| \leq  \frac{2D |\nabla u|}{u} D^2 \abs{\nabla \ln\br{\frac{u}{D}}}^2 + 2D \abs{\nabla \ln\br{\frac{u}{D}}} D^2 \abs{\nabla^2 \ln\br{\frac{u}{D}}} \\
 \lesssim D \abs{\nabla \ln\br{\frac{u}{D}}} + D^2 \abs{\nabla^2 \ln\br{\frac{u}{D}}}
\end{multline*}
by using once \eqref{cond.uD} and twice \eqref{eq.DlnuD}. The bounds on $T_2$ and $T_3$ are equally easy: after applying the chain rule, we will have a sum of terms that will contain either $D \abs{\nabla \ln\br{\frac{u}{D}}}$,  $D^2 \abs{\nabla^2 \ln\br{\frac{u}{D}}}$ or $D\abs{\nabla(\abs{\nabla D})}$, and we bound the rest using \eqref{cond.uD} or \eqref{eq.DlnuD}. \ep

\medskip

\noindent {\em Proof of Theorem \ref{Maincod1} $(i)\implies(ii)$: } Let $\beta>0$ be fixed. When $d=n-1$, we showed in Corollary \ref{CorDb2} that $D_\beta\abs{\nabla^2D_\beta}$ satisfies the Carleson measure condition in $\om$. Then the implication $(i)\implies(ii)$ of Theorem \ref{Maincod1} is a straightforward consequence of \eqref{eq.D2u}, Theorem \ref{ThNlnu/D}, Corollary \ref{CorCaccio}, and \eqref{Db=dist}. \ep

\medskip
In the general case when $0<d<n$, $D_\beta\abs{\nabla^2D_\beta}$ is no longer the right quantity to look at. Instead, we know that $D_\beta\abs{\nabla(\abs{\nabla D_\beta}^2)}\in CM_\om$ by Theorem \ref{thmDEM} or Corollary \ref{CorDb3}.
\medskip

\noindent {\em Proof of Theorem \ref{mainthm} $(i)\implies(iii)$: } An immediate consequence of Theorem \ref{thmDEM}, Theorem \ref{ThNlnu/D}, Corollary \ref{CorCaccio}, and the inequality \eqref{eq.D|Du|2}.
 \ep

\medskip

We now show that the Carleson measure estimate can take a simpler form, namely, we show that $(i)\implies(ii)$ in Theorem \ref{mainthm} holds. We aim to use the relation
\begin{equation}\label{eq.DDu=DDu2}
     \frac{D^2\abs{\nabla(\abs{\nabla u})}}{u}=\frac{u}{2D\abs{\nabla u}}\cdot\frac{D^3\abs{\nabla\br{\abs{\nabla u}^2}}}{u^2},
\end{equation}
and the fact that $(i)\implies(iii)$ in Theorem \ref{mainthm} holds.
The issue with this strategy is that since $\abs{\nabla u}$ might be 0 in $\om$, we cannot simply pull out $\frac{u}{2D\abs{\nabla u}}$ as a constant when taking integrals on both sides of \eqref{eq.DDu=DDu2}. But fortunately, we have the following lemma which says that the set where $\frac{u}{D\abs{\nabla u}}$ is big is small in the sense that it satisfies a Carleson condition.

\begin{lemma}\label{lem.1u/DDuCM}
       Let $d<n$, and let $\om$ be a domain in $\Rn$ with uniformly $d$-rectifiable boundary $\pom$. Let $L=-\diver \br{D_\beta^{d+n-1}\A\nabla}$ be uniformly elliptic. 
Assume either
\begin{enumerate}
    \item $d=n-1$, $\Omega$ is uniform, and $\A$ satisfies the DKP condition, or 
    \item $0<d<n-1$ and $\A=I$.
\end{enumerate}  
Then for any $c'\in(0, c/2]$, where $c$ is the constant in Corollary \ref{CorDb3}, there exists $C>0$ such that for any $x\in\pom$, any $r\in(0,\diam\pom)$, and any positive solution $u$ to $Lu=0$ in $B(x,2r)$ with $u=0$ on $\pom\cap B(x,2r)$, we have 
    \[
\iint_{B(x,r)\cap\om}\1_{\set{\frac{D_\beta\abs{\nabla u}}{u}<c'}}(X)\delta_{\pom}(X)^{d-n}dX\le C\,r^d.
    \]
\end{lemma}

\bp
For simplicity, let us write $D$ for $D_\beta$. For any $\eta>0$, denote 
\[E_\eta:=\set{\abs{\frac{D\nabla u}{u}-\nabla D}>\eta}=\set{D\abs{\nabla \ln\br{\frac{u}{D}}}>\eta}.\]
Then for any $\eta>0$, 
\begin{equation}\label{1EKCM}
    \eta^2\iint_{B(x,r)\cap\om}\1_{E_\eta}(X)\delta_{\pom}(X)^{d-n}dX
    <\iint_{B(x,r)\cap\om}D^2\abs{\nabla \ln\br{\frac{u}{D}}}^2\delta_{\pom}^{d-n}dX
    \le C\,r^d,
\end{equation}
    where we have used Theorem \ref{ThNlnu/D} and \eqref{Db=dist} in the last inequality. Let $c$ be the constant in Corollary \ref{CorDb3}. By the triangle inequality, we have 
\[
\set{\frac{D\abs{\nabla u}}{u} < \frac{c}{2}}\setminus\set{\abs{\nabla D}<c}\subset \set{\abs{\frac{D\nabla u}{u}-\nabla D}>\frac{c}{2}}=E_{\frac{c}{2}},
\]
which implies that 
\[\1_{\set{\frac{D\abs{\nabla u}}{u}< \frac{c}{2}}}\le \1_{\set{\abs{\nabla D}<c}}+\1_{E_{\frac{c}{2}}}.\]
Thus, the lemma follows from \eqref{1EKCM} and Corollary \ref{CorDb3}.   
\ep

\medskip

Now we are ready for \smallskip

\noindent {\em Proof of Theorem \ref{mainthm} $(i)\implies(ii)$: }
Let $K>0$ be a (large) constant to be determined later. We write 
\begin{multline*}
    \iint_{B(x,r)\cap\om} \frac{\abs{\nabla\abs{\nabla u}}^2}{u^2}\delta_{\pom}^{d+4-n}dX
=\iint_{B(x,r)\cap \set{\frac{u}{D\abs{\nabla u}}>K}} \frac{\abs{\nabla\abs{\nabla u}}^2}{u^2}\delta_{\pom}^{d+4-n}dX \\
+\iint_{B(x,r)\cap \set{\frac{u}{D\abs{\nabla u}}\le K}} \frac{\abs{\nabla\abs{\nabla u}}^2}{u^2}\delta_{\pom}^{d+4-n}dX =:I_1+I_2,
\end{multline*}
where we write $D$ for $D_\beta$ for simplicity. 
For $I_2$, we use \eqref{eq.DDu=DDu2} to get that
\[
    I_2\le CK^2 \iint_{B(x,r)\cap\om}\frac{\abs{\nabla\br{\abs{\nabla u}^2}}^2}{u^4}\delta_{\pom}^{d+6-n}dX\le CK^2r^d
\]
by $(i)\implies(iii)$ in Theorem \ref{mainthm}. To treat $I_1$, we first observe that the set 
\[
S_K:=B(x,r)\cap\set{\frac{u}{D\abs{\nabla u}}>K}=B(x,r)\cap\set{\frac{D\abs{\nabla u}}{u}<\frac{1}{K}}
\]
is open as $\frac{D\abs{\nabla u}}{u}$ is (H\"older) continuous (see Lemma \ref{lemDNu/uHolder}). Thus we can cover $S_K$ by Whitney balls:
\[S_K\subset \bigcup_{X\in S_K}B(X,c_0\delta_{\pom}(X)),\]
 where $c_0<1/2$ is some small constant to be determined. By the Vitali covering lemma, there exists a family of disjoint balls $B_i:=B(X_i,c_0\delta_{\pom}(X_i))$, $i\in \mathcal{I}$, such that $S_K\subset\bigcup_{i\in\mathcal{I}}5B_i$. Then by Lemma \ref{lemuinW22}, the Harnack principle, Lemma \ref{lemDNu/u<1}, and disjointness of $\set{B_i}_{i\in\mathcal{I}}$, we have 
 \begin{multline*}
     I_1\le \sum_{i\in\mathcal{I}}\iint_{5B_i}\frac{\abs{\nabla\abs{\nabla u}}^2}{u^2}\delta_{\pom}^{d+4-n}dX
     \le C\sum_{i\in\mathcal{I}}\iint_{10B_i}\frac{\abs{\nabla u}^2}{u^2}\delta_{\pom}^{d+2-n}dX\le C\sum_{i\in\mathcal{I}}\iint_{10B_i}\delta_{\pom}^{d-n}dX\\
     \le C\sum_{i\in\mathcal{I}}\delta(X_i)^d
     \le C\sum_{i\in\mathcal{I}}\iint_{B_i}\delta_{\pom}^{d-n}dX\le C\iint_{\bigcup_{i\in\mathcal{I}}B_i}\delta_{\pom}^{d-n}dX.
 \end{multline*}
Thanks to Lemma \ref{lemDNu/uHolder}, for any $K>0$, we can choose $c_0$ sufficiently small so that 
\[
\frac{D(X)\abs{\nabla u(X)}}{u(X)}<\frac{2}{K}\quad \text{for any }X\in \bigcup_{i\in\mathcal{I}}B_i.
\]
Also, we have $\bigcup_{i\in\mathcal{I}}B_i\subset B(x,2r)$. Take $K>4/c$, where $c$ is the constant in Corollary \ref{CorDb}. Then 
\[
I_1\le C\iint_{B(x,2r)}\1_{\set{\frac{D\abs{\nabla u}}{u}<\frac{c}{2}}}\delta_{\pom}^{d-n}dX\le Cr^d
\]
by Lemma \ref{lem.1u/DDuCM}. This completes the proof.
\ep

\section{The limit case}\label{SecLimit}

We want to prove that good estimates on the Green functions imply uniform rectifiability of the boundary. We shall use a compactness argument (similar, for instance, to \cite{DEM}, \cite{DMgreen}, \cite{FLM}), that is, we shall ``blow-up'' some solutions and reduce matters to the limit case, which is to show that if $|\nabla G|$ is constant, then the boundary is flat - i.e. up to translations or dilatations, the domain is either $\R^n_+$, $\{(x,t) \in \R^{n-1} \times \R, \, |t|>1\}$, or $\R^n \setminus \R^d$ for some integer $d$.

\begin{theorem}[Theorem 3.1 in \cite{DEM}] \label{ThDEM}
Let $E \subset \R^n$ be a closed set, and let $G$ be a continuous non-negative function on $\R^n$ that vanishes on $E$, is of class $C^1$ on $\Omega := \R^n \setminus E$, and such that $|\nabla G|$ is positive and constant on every connected component of $\Omega$. Then the set $E$ is convex. If, moreover, $|\nabla G|=1$ on $\Omega$, then $G(X) = \dist(X,E)$ for $X\in \R^n$.
\end{theorem}

Let us give two intermediate propositions. The first is also from \cite{DEM} and treats the case where $\Omega = \R^n \setminus \partial \Omega$:

\begin{proposition} \label{propDEM}
Let $d<n$ and $E\subset \R^n$ be a convex, $d$-Ahlfors regular set, such that $\delta_E(X) = \dist(X,E)$ is of class $C^2$ on $\Omega := \R^n \setminus E$. Then $d$ is an integer and $E$ is a $d$-dimensional plane.
\end{proposition}

\bp A proof of this result is given in the proof of \cite[Corollary 3.2]{DEM}. Here, we present a different proof, which is more direct to the authors. 

If $E$ is $d$-Ahlfors regular, and $m$ is an integer such that $d>m$, then we can find $m+1$ points $\{x_i\}_{1\leq i \leq m+1}$ in $E$ that are not contained in a $m$-plane. But since $E$ is convex, $E$ has to contain the convex hull of the $\{x_i\}_i$, hence has to be at least of dimension $m+1$. We deduce that $d$ is an integer and $E$ is included in a $d$-plane $P_E$.

It remains to prove that $E = P_E$. Assume by contraposition that it is not the case, and we shall prove that $\delta_E$ is not of class $C^2$ in $\Omega$. Without loss of generality, we can assume $0 \in E$, $P_E =  \R^{d} \times \{0_{\R^{n-d}}\}$, and $E \subset \R_+ \times \R^{d-1} \times \{0_{\R^{n-d}}\}$ (note that we necessarily have $d\ge1$ otherwise $E=\{0\}=P_E$). We pretend that we can also assume that $e_1 := (1,0,\dots,0) \in E$. Indeed, if it is not the case, we dilate and rotate $E$ with the following algorithm:
\begin{enumerate}
\item If there is a $r$ such that $(r,0,\dots,0) \in E$ and so - up to dilatation - $e_1 \in E$.
\item If $E \subset (\R_+)^j \times \R^{d-j} \times \{0_{\R^{n-d}}\}$ and $E \cap (\R^{j} \times \{0_{\R^{n-j}}\}) = \{0\}$. By the Hahn-Banach theorem and up to a rotation, 
\begin{equation} \label{condonErot}
E \subset (\R_+)^{j+1} \times \R^{d-j-1} \times \{0_{\R^{n-d}}\}.
\end{equation} 
If $E \cap (\R^{j+1} \times \{0_{\R^{n-j-1}}\}) \supset \{0,e\}$ with $e \neq 0$, we rotate $E$ around $\{0_{\R^{j+1}}\} \times \R^{n-j-1}$ and then dilate so that $e = e_1$, and we observe that \eqref{condonErot} forces $E$ to stay in $\R_+ \times \R^{d-1} \times \{0_{\R^{n-d}}\}$ after rotation. Otherwise, we repeat (2) where $j$ becomes $j+1$.
\item We have $E \cap (\R^{d} \times \{0_{\R^{n-d}}\}) = E \neq \{0\}$, so the algorithm necessarily stops looping when $j=d-1$.   
\end{enumerate}
With all this, if $e_{n} = (0,\dots,0,1)$ is a vector normal to $P_E$, we can check that 
\[\delta_E(te_1 + e_{n}) = |(te_1 + e_{n}) - (te_1)| = 1 \qquad \text{ when } t\in [0,1]\]
and
\[\delta_E(te_1 + e_{n}) = |(te_1 + e_{n}) - 0| = \sqrt{1+t^2} \qquad \text{ when } t\leq0.\]
Differentiating twice in $t$ shows that $t\mapsto \delta_E(te_1 + e_{n}) $ is not $C^2$ at $0$, therefore $\delta_E$ is not of class $C^2$ around $te_1 + e_{n}$.
\ep

If $\Omega \subsetneq \R^n \setminus \partial \Omega$ - which means that $d\geq n-1$ by the proof of Proposition \ref{lowdim=>uniform} - we shall use the following result instead.

\begin{proposition} \label{propDEMb}
Let $E \subset \R^n$ be a closed convex set with $(n-1)$-dimensional $C^3$-boundaries (i.e. 
$\partial E$ can be locally described as the graph of a $C^3$ function with $n-1$ variables). Write $\Omega$ for $\R^n \setminus E$ and $\delta_E$ for the function $\dist(.,E)$. If there exists an elliptic operator $L=-\div \A \nabla$ {\bf with constant coefficients} such that $L \delta_E = 0$ on $\Omega$, then each connected component of $\Omega$ is $\R^n_+$, up to rotation and translation. 
\end{proposition}

\bp Take a point $X \in \Omega$. We write $x \in \partial E$ for its projection on $E$, and $t$ for its distance, i.e. $|X-x| = t = \delta_E(X)$. We parametrize the boundary $\partial E$ around $x$ by a concave function $\varphi \in C^3$, that is, for a small $r$, we have 
\[E \cap B(x,r) = \{(y,s) \in \R^{n-1} \times \R, \, s\leq \varphi(y)\} \cap B(x,r).\] Up to rotation, we can always assume that 
$\nabla \varphi(x) = 0$ - that is, $X$ is ``above'' $x$ - and the Hessian matrix $D^2\varphi(x)$ is diagonal, with eigenvalues $-\lambda_i \leq 0$ on the diagonal. The proposition will be proved if we show that
\begin{equation*}
L \delta_E(X) = 0 \implies \lambda_i = 0 \text{ for all } i\in \{1,\dots,n-1\}.
\end{equation*}

We give a proof without subtleties. Since $\varphi$ is convex and $C^2$, the map
\[
\rho(y,s):= (y,\varphi(y)) + s \frac{(-\nabla \varphi(y),1)}{|(-\nabla \varphi(y),1)|}.
\]
is a bijection from a neighborhood of $(x,t)$ to a neighborhood of $X$. Moreover, in this neighborhood, we have $\rho^{-1}(Y) = (y,\delta_E(Y))$, where $y$ is such that $(y,\varphi(y))$ is the (unique) projection of $Y$ on $E$. As a consequence, the function $\delta_E$ is given around $X$ by $\pi_n \circ \rho^{-1}$, where $\pi_n:\, (y,s) \in \R^{n-1} \times \R \mapsto s \in \R$ is the projection on the last coordinate. We shall brutally compute $L\delta_E(X)$ (note that assuming $\nabla \varphi(x) = 0$ greatly simplifies the computations). We start with computing the Jacobian matrix of $\rho$ at $(x,t)$: 
\begin{equation} \label{defJacZ}
\Jac_\rho(x,t) = \nabla \rho(x,t)=
\begin{bmatrix}
\begin{BMAT}{c.c}{c.c}
I_{n-1} & \mathbf{0}  \\
\mathbf{0} & 1
\end{BMAT}
\end{bmatrix}+
\begin{bmatrix}
\begin{BMAT}{c.c}{c.c}
-t  D^2\varphi(x) & \mathbf{0}  \\
\mathbf{0} & 0
\end{BMAT}
\end{bmatrix},
\end{equation}
that is, $\Jac_\rho(x,t)$ is a diagonal matrix with entries $(1+t\lambda_1),\dots,(1+t\lambda_{n-1}),1$. Since the unit vector $e_n = (0,\dots,0,1)$ verifies
$e_n = \nabla\pi_n=\nabla [\delta_E \circ \rho] = \Jac_\rho [\nabla \delta_E] \circ \rho$, 
we deduce
\begin{equation} \label{NdeltaEZ}
\nabla \delta_E (X) = [\nabla \delta_E]\circ \rho(x,t) = \Jac_\rho(x,t)^{-1} e_n = e_n.
\end{equation}
Next, we compute the second derivatives of $\delta_E$. To this end, we take the second derivatives of $\rho$, which is valid as $\vp\in C^3$. A direct computation yields
\[ 
\left\{ \begin{array}{ll}  
\partial_i \partial_j \rho(x,t) = \big(-t\partial_i \partial_j \nabla \varphi(x), \partial_i\partial_j \varphi(x)\big) & \text{ if } 1 \leq i,j \leq n-1 \\
\partial_i \partial_n \rho(x,t) = \big(-\partial_i \nabla \varphi(x), 0\big) & \text{ if } 1 \leq i \leq n-1 \\
\partial_n\partial_n \rho (x,t) = 0.
 \end{array} \right.
 \]
 Applying the chain rule, we have
\begin{multline}\label{eq.DEMb1}
0 = \partial_i \partial_j\pi_n=\partial_i \partial_j[\delta_E \circ \rho] = \partial_i \{(\partial_j \rho) \cdot ([\nabla \delta_E] \circ \rho)\} \\
= (\partial_i\partial_j \rho) \cdot ([\nabla \delta_E] \circ \rho) + (\partial_j \rho) \br{[D^2 \delta_E] \circ \rho} \cdot (\partial_i \rho).
\end{multline}
Set $\lambda_n=0$, then for $i,j\in\set{1,\dots,n}$,
\[
(\partial_j \rho)(x,t) [D^2 \delta_E] \circ \rho(x,t) \cdot (\partial_i \rho)(x,t)=(1+t\lambda_i)(1+t\lambda_j) [\partial_i \partial_j \delta_E](X).
\]
Hence, by \eqref{eq.DEMb1}, \eqref{NdeltaEZ}, and the computations of $\partial_i\partial_j\rho(x,t)$, we obtain
\begin{multline*}
(1+t\lambda_i)(1+t\lambda_j) [\partial_i \partial_j \delta_E](X) = - \{\partial_i\partial_j \rho(x,t)\} \cdot \{[\nabla \delta_E](X)\}\\
= -\partial_i \partial_j \varphi(x)  = \left\{ \begin{array}{ll} \lambda_i = \lambda_j & \text{ if } i = j < n\\ 0 & \text{ if } i \neq j \text{ or } i =j= n. \end{array} \right.
\end{multline*}
We just showed that 
\[(\partial_i \partial_j \delta_E)(X) = \frac{\delta_{ij} \lambda_i}{(1+t\lambda_i)^2} \geq 0.\]
If $a_{ij}$ are the entries of constant matrix $\A$, then $a_{ii}>0$ by ellipticity, and thus 
\[L\delta_E(X) = - \sum_{i=1}^{n-1} a_{ii} \frac{\lambda_i}{(1+t\lambda_i)^2} \leq 0. \]
Equality holds if and only if $\lambda_i = 0$ for all $i\in \{1,\dots,n-1\}$. The  proposition follows.
\ep

We invoke the last three results to obtain the following corollary.

\begin{corollary} \label{corDEM} Let $d< n$ ($d$ is not necessarily an integer), $E \subset \R^n$ be a closed set with $d$-Ahlfors regular boundaries. Let $\A$ be a positive definite constant matrix, and $D\in C^\infty(\R^n\setminus E)$ be such that 
\[C^{-1} \dist(X,E) \leq D(X) \leq C \dist(X,E) \quad \text{ for } X\in \Omega := \R^n \setminus E\]
and $L = - \div\br{D^{d+1-n}\A \nabla}$. Assume the existence of a positive solution $G$ to $LG=0$ on $\Omega$ that continuously vanishes on $E$ and that $|\nabla G|$ is constant on every connected component of $\Omega$. Then $d$ is an integer. Moreover, if $d<n-1$, then $E$ is a $d$-dimensional plane, and if $d=n-1$, then $E$ is either a half space, a $(n-1)$-plane, or the set between two parallel $(n-1)$-planes.
\end{corollary}

\bp
Since $G$ is a solution of a uniformly elliptic operator with smooth (even constant when $d=n-1$) coefficients , $G \in C^\infty(\Omega)$. Moreover, the constant $|\nabla G|$ is necessarily positive, since $G$ is positive on $\Omega$ and is zero on $\partial E$. Theorem \ref{ThDEM} yields that $E$ is convex and - up to a multiplicative constant depending only on the connected component - we have $G(X) = \dist(X,E)$. 

In the case where $E$ has empty interior - meaning that $E = \partial E$ - we have that $E$ is $d$-Ahlfors regular, and $\dist(X,E) = G(X) \in C^\infty(\Omega) \subset C^2(\Omega)$, so Proposition \ref{propDEM} applies. We deduce that $d$ is an integer and $E$ is a $d$-plane.

If $E$ has non-empty interior, then $E$ has dimension $n$. Also, $\partial E$ is locally Lipschitz as $E$ is convex. It entails that  $d=n-1$, and thus $L$ has constant coefficients. For $t>0$, we define $E_t = \{X\in \R^n, \, \dist(X,E) \leq t\}$. Note that $E_t$ is a set with $C^\infty$ boundary since the boundary is the level set of $G$, which is a smooth function whose gradient is never 0 on $\om$. Define 
\[
G'(X)=
\begin{cases}
G(X)-t &\quad \text{if }X\in \om_t:=\Rn\setminus E_t\\
0 &\quad \text{if } X\in E_t.
\end{cases}
\]
Observe that $G'$ is a nonnegative and continuous function in $\Rn$ which vanishes on $E_t$, $G'\in C^\infty(\om_t)$, and $\abs{\nabla G'}$ is constant on every component of $\om_t$ because $\nabla G'=\nabla G$ on $\om_t$. So by Theorem \ref{ThDEM}, $G' = c\dist(.,E_t) =: c \delta_{E_t}$ on each connected component of $\Omega_t$. Since $LG' = LG=0$ on $\om_t$, we have $L\delta_{E_t} = 0$ on $\Omega_t$. Proposition \ref{propDEMb} gives that each connected component of $\Omega_t$ is $\R^n_+$ up to rotation and translation. Either $\Omega_t$ is connected and $\Omega_t$ is $\R^n_+$. Or else $\Omega_t$ is disconnected, and it has exactly two connected components which are two copies of $\R^n_+$; since the two copies are non-overlapping, it means that the boundaries of $\Omega_t$ are parallel $(n-1)$-planes, and $E_t$ is the set between those two planes. 

Once you know $E_t$ for one $t>0$, you also know $E$: if $E_t$ is a half space, then $E$ is a half space obtained by translating $E_t$ by $t$ units in the direction of $-\nabla G$; and if $E_t$ is between two parallel planes $P_{1,t}$ and $P_{2,t}$, then $E$ is the space between $P_1$ and $P_2$, where $P_i$ is obtained by translating $P_{i,t}$ by $t$ units in the direction of $-\nabla G(X_{i})$, with $X_i \in P_{i,t}$. 
\ep

\section{Control on the second derivatives on the Green function implies UR}\label{SecGtoUR}

Let us introduce the local Hausdorff distance for sets as 
\begin{multline*}
d_{x,r}(E,F) := \frac1r \sup\{ \dist(y,F), \, y \in E \cap B(x,r)\} + \frac1r \sup\{ \dist(y,E), \, y \in F \cap B(x,r)\}.
\end{multline*}

\subsection{The compactness argument}

\medskip
We want to prove the convergence of a sequence of domains $\om_j$, a sequence of operators $L_j=-\diver(D_{\beta,j}^{d+1-n}\A_j\nabla )$, and the corresponding sequence of solutions $u_j$, and to show that their limits have the desired properties. Some of the arguments in this subsection are similar to the ones in \cite[Theorem 4.8]{HMMTZ} or \cite[Theorem 2.19]{DMgreen}.

\begin{proposition} \label{problowup1}
Let $d<n$ and let $\{\Omega_j\}_{j\geq 0}$ be a sequence of domains in $\R^n$ such that
\begin{enumerate}[(i)]
\item $\Omega_j$ is uniform with (uniform) constants $(\epsilon,N)$;
\item the boundary $\partial \Omega_j$ is $d$-Ahlfors regular, and the $d$-Ahlfors regular measure $\sigma_j$ satisfies \eqref{defADR} with a uniform constant $C$;
\item  $\diam \partial \Omega_j \geq 2^j$;
\item $0\in \partial \Omega_j$.
\end{enumerate}
Then, up to a subsequence, $\Omega_j$ converges to a domain $\Omega_\infty$, in the sense that, for any $x\in \R^n$ and $r>0$,
\[\lim_{j\to \infty} d_{x,r}(\Omega_j,\Omega_\infty) \to 0 \quad \text{ and } \quad \lim_{j\to \infty} d_{x,r}(\partial \Omega_j,\partial \Omega_\infty) \to 0.\]
The limit domain $\Omega_\infty$ is uniform with constants $(\epsilon/4,10N)$. The measure $\sigma_j$ converges weakly to a $d$-Ahlfors regular measure $\sigma_\infty$ whose support is unbounded and equal to $\partial \Omega_\infty$, meaning in particular that $\partial \Omega$ is $d$-Ahlfors regular. For $\beta >0, \, j \in \mathbb N \cup \{\infty\}$, define the smooth distance $D_{\beta,j}$ as
\[
D_{\beta,j}(X) := \left( \int_{\partial \Omega_j} |X-y|^{-d-\beta} d\sigma_j(y) \right)^{-\frac{1}{\beta}}
\]
for $X \in \R^n \setminus \partial \Omega_j$, and $D_{\beta,j}(X) = 0$ for $X\in \partial \Omega_j$.
Then $D_{\beta,j}$ converges uniformly to $D_{\beta,\infty}$ on any compact  sets of $\R^n$.
\end{proposition}

\bp
By the Blaschke selection theorem (see \cite[Theorem 3.16]{Falconer}), there exists a subsequence - still called $\Omega_j$ - and two closed sets $E$ and $F$ such that, for $x\in \R^n$ and $r>0$, 
\[\lim_{j\to \infty} d_{x,r}(\R^n \setminus \Omega_j,E) \to 0 \quad \text{ and } \quad \lim_{j\to \infty} d_{x,r}(\partial \Omega_j,F) \to 0.\]
Define $\om_\infty:=\Rn\setminus E$; we want to prove that $\partial E = F$. We easily have $\partial E \subset F$, since $z\in \partial E$ means that we can find values of $\Omega_j$ as close as we want to $z$. 
To see that $F\subset\partial E$, fix any  $z\in F$ and $r>0$, we show that $B(z,r)\cap\om_\infty\neq\emptyset$. We take a sequence $z_j \in \partial \Omega_j$ such that $z_j \to z$, and then a corkscrew point $a_j$ in $\Omega_j \cap B(z_j,3r/4)$; for $j$ large enough, $a_j \in B(z,r)$ and $B_j:= B(a_j,\epsilon r/4) \subset \Omega_\infty$, and any of such $a_j$ is a corkscrew point for $\Omega_\infty$. This proves that $F \subset \partial E$ and $\Omega_\infty$ satisfies the corkscrew point condition with a constant $\epsilon/4$. 

The Harnack chain is also simple. We write $\delta_j$ for $\delta_{\Omega_j} = \dist(.,\partial \Omega_j)$. If $X,Y\in \Omega_\infty$ are such that $|X-Y|/\min\{\delta_\infty(X),\delta_\infty(Y)\} \leq 40/\epsilon$, then we also have $|X-Y|/\min\{\delta_j(X),\delta_j(Y)\} \leq 50/\epsilon$ for $j$ large enough. Those can be linked by an Harnack chain $\{X_i\}_{0\leq i \leq 5N}$ inside $\Omega_j$  (see Proposition \ref{propHarnack}). We know that $|X_i - X_{i+1}| \leq \delta_j(X_i)/2$: either it is also true with $\delta_\infty$ too, or we add the middle point of $[X_i,X_{i+1}]$ in the chain. 
We refer a reader interested in more details for the argument given so far to the proof of \cite[Theorem 4.8 (v)]{HMMTZ}.

Now, we want to show that $\partial \Omega_\infty$ is unbounded and $d$-Ahlfors regular. The first fact is immediate, since 
\[\diam \partial \Omega_\infty = \diam F = \lim_{j\to \infty} \diam \partial \Omega_j = \infty.\]
As for the second fact, up to a subsequence, $\sigma_j$ converges weakly-$*$ to a Radon measure $\sigma_\infty$. The fact that $\sigma_\infty$ is a $d$-Ahlfors regular measure with support on $\partial \Omega_\infty$ is 
\cite[Theorem 4.8 (v)]{HMMTZ}.

It remains to prove the uniform convergence of $D_{\beta,j}$ to $D_{\beta,\infty}$ on each ball $B_k := B(0,2^k)$. Since the $d$-Ahlfors regular constants are uniform, there is a $K = K(\beta,C)\ge1$ such that,
\begin{equation} \label{Db=dist2}
K^{-1}\delta_j(X) \leq D_{\beta,j}(X) \leq K\delta_j(X) \qquad \text{ for } X\in \R^n, \, j\in \mathbb N \cup \{\infty\},
\end{equation}
see \eqref{Db=dist}. Let $\epsilon>0$ and take $\delta:= \epsilon/5K$. Then, for $j$ large enough,
\[|D_{\beta,j}(X) - D_{\beta,\infty}(X)| \leq \epsilon \qquad \text{ for $X\in \R^n$ such that } \delta_\infty(X) \leq \delta.\] 
For $X\in B_k$ and $m\geq 2k$, by decomposing $\pom_j\setminus B_m=\bigcup_{i=0}^\infty\pom\cap \br{B_{m+i+1}\setminus B_{m+i}}$, one can show that
\begin{equation}\label{eq.Dbetajm1}
    \int_{\partial \Omega_j \setminus B_m} |X-y|^{-d-\beta} d\sigma_j(y) \leq K' 2^{-\beta m}
\end{equation}
for some $K'=K'(\beta,C,d)>0$. Define
\[D_{\beta,j,m}(X) := \left( \int_{\partial \Omega_j \cap B_m} |X-y|^{-d-\beta} d\sigma_j(y) \right)^{-\frac1\beta}. \]
By considering the function $z\mapsto z^{-1/\beta}$, we have 
\begin{multline*}
    \abs{D_{\beta,j,m}(X) - D_{\beta,j}(X)}
    \leq\frac{1}{\beta}\br{D_{\beta,j,m}^{-\beta}(X)}^{-\frac{1}{\beta}-1} \abs{D_{\beta,j,m}^{-\beta}(X) - D_{\beta,j}^{-\beta}(X)}\\
    \leq C_{d,\beta} 2^{k(\beta+1)}\int_{\partial \Omega_j \setminus B_m} |X-y|^{-d-\beta} d\sigma_j(y)\le C_{d,\beta,K'}2^{k(1+\beta)}2^{-\beta m}
\end{multline*}
due to \eqref{eq.Dbetajm1} and 
$D_{\beta,j,m}^{-\beta}(X)\ge \int_{\pom_j\cap B_k}|X-y|^{-d-\beta} d\sigma_j(y)\ge 2^{-d-\beta}2^{-k\beta}$.
This implies that we can picking $m>k(1+\beta)/\beta$  sufficiently large so that
\[|D_{\beta,j,m}(X) - D_{\beta,j}(X)| + |D_{\beta,\infty,m}(X) - D_{\beta,\infty}(X)| \leq \epsilon/2 \qquad \text{ for } X\in B_k.\]

Now, we conclude that, for $X\in B_k$ with $\delta_\infty(X) \geq \delta$,
\begin{multline*}
|D_{\beta,j}(X) - D_{\beta,\infty}(X)| \leq \frac12 \epsilon + |D_{\beta,j,m}(X) - D_{\beta,\infty,m}(X)|   \\
 \leq \frac12 \epsilon + C_{d,\beta}2^{k(\beta+1)} \left| \int_{\partial \Omega_j \cap B_m} |X-y|^{-d-\beta} d\sigma_j(y) - \int_{\partial \Omega_\infty \cap B_m} |X-y|^{-d-\beta} d\sigma_\infty(y) \right|  \leq \epsilon
\end{multline*}
if $j$ is large enough, since $\sigma_j$ converges weakly-$*$ to $\sigma_\infty$ and the functions $|X-y|^{-d-\beta}$ are uniformly bounded on a neighborhood of $\partial \Omega_\infty$. This proves the uniform convergence of $D_{\beta,j}$ to $D_{\beta,\infty}$ on $B_k$ and hence concludes the proposition.
\ep

\begin{proposition} \label{problowup2}
Let $\beta>0$, $\{\Omega_j\}_{j\geq 1}$ and $D_{\beta, j}$ be as in Proposition \ref{problowup1}. Define $\delta_j$ as the distance function to $\partial \Omega_j$, $B_j$ for $B(0,2^j)$, and $\Omega_{j,j}:= \{X \in \Omega_j \cap B_j, \, \delta_j(X) \geq 2^{-j}\}$.

\medskip

Let $L_j:= -\div [D_{\beta,j}^{d+1-n} \A_j \nabla]$ be a sequence of uniformly elliptic operators\footnote{We mean that the ellipticity and boundedness constants of $\A_j$ are independent of $j$ too.} on $\Omega_j \cap B_j$ that satisfies 
\begin{equation} \label{assumptionblowup3}
\sup_{j\in \mathbb N} \|\delta_j \nabla \A_j\|_{L^\infty(\Omega_j)} < \infty
\end{equation}
and
\begin{equation} \label{assumptionblowup2}
\iint_{\Omega_{j,j}}|\nabla \A_j| dX \leq 2^{-j}.
\end{equation}
Assume in addition that there is a positive weak solution $u_j$ to $L_j u = 0$ in $\Omega_j \cap B_j$ that satisfies $u_j = 0$ on $\partial \Omega_j \cap B_j$ and 
\begin{equation} \label{assumptionblowup}
\iint_{\Omega_{j,j}} \frac{\abs{\nabla \br{|\nabla u_j|^2}}^2}{u_j^4} dX \leq 4^{-j}.
\end{equation}

Then, up to a subsequence, $L_j$ converges to an uniformly elliptic operator $L_\infty$ on $\Omega_\infty$ that satisfies $ L_\infty := - \div [D_{\beta,\infty}^{d+1-n} \A_\infty \nabla]$ for a constant matrix $\A_\infty$. In addition, there exists a positive continuous solution $u_\infty$ to $L_\infty u = 0$ on $\Omega_\infty$ such that $u=0$ on $\partial \Omega_\infty$ and $|\nabla u| = 1$ on $\Omega_\infty$.
\end{proposition}

\bp
Define $\Omega_{\infty,j}$ as $\{X \in \Omega_\infty \cap B_{j-1}, \, \delta_j(X) \geq 2^{1-j}\}$. Without loss of generality (i.e. up to taking a subsequence of the $\Omega_j$), we can assume that $\Omega_{\infty,j} \subset \Omega_{j,j}$.

By \eqref{assumptionblowup3}, for any $k\geq 1$, the matrices $\A_j$ are equicontinuous on $\Omega_{\infty,k}$. Together with the uniform uniform ellipticity and the Arzel\`a-Ascoli theorem, we deduce that, up to a subsequence, the matrices $\A_j$ converge to a matrix $\A_{k,\infty}$ on $\Omega_{\infty,k}$. A diagonal process gives that $\A_j$ converges to a matrix $\A_\infty$ uniformly on compact sets of $\Omega_\infty$. By taking the extra subsequence from Proposition \ref{problowup1}, we obtain that 
\begin{equation} \label{AjtoAinfty}
    A_j := D^{d+1-n}_{\beta,j} \A_j \text{ converges to } A_\infty := D^{d+1-n}_{\beta,\infty} \A_\infty \text{ uniformly on compact sets of $\Omega_\infty$.} 
\end{equation}

Let $X_0 \in \Omega_\infty$ be a corkscrew point for the ball $B(0,1)$. Without loss of generality, we can extend $u_j$ to $B_j$ by setting $u_j = 0$ on $\R^n \setminus \Omega_j$ and we can rescale $u_j$ so that $u_j(X_0) = 1$. Since the $u_j$ are H\"older continuous on compact sets of $\R^n \cap B_{j-1}$\footnote{The H\"older regularity constants depend on $u_j$ only via $u_j(X_0)$ - since the $\Omega_j$ are uniformly uniform and the $L_j$ are uniformly elliptic - so there are independent of $j$. See for instance Lemmas 11.30, 11.32, and 15.14 in \cite{DFMmixed}.}, the Arzel\`a-Ascoli theorem implies that, up to a subsequence, $u_j$ converges uniformly on any compact set of $\R^n$ to a function $u_\infty$; in particular, $u_\infty = 0$ on $\partial \Omega_\infty$.
We claim that $u_j$ also converges to $u_\infty$ in  $W^{1,2}_{loc}(\Omega_\infty)$. The $L^2$ convergence follows from uniform convergence, and so we only need to show that for any compact subset $K\Subset\om_\infty$,  
\begin{equation}\label{eq.localW12conv}
    \iint_K\abs{\nabla(u_j-u_k)}^2dX \to 0 \text{ as }j,k\to\infty.
\end{equation}
Let $\vp\in C_0^\infty(\om_\infty)$ be a nonnegative function such that $\vp=1$ on $K$. By the ellipticity of $\A_j$ and the fact that $D_{\beta,j}^{d+1-n}\approx\delta_j^{d+1-n}\approx  C_K$ for some positive constant $C_K$ when $j$ is sufficiently large, we have \begin{multline*}
 \iint_{\om_\infty}\abs{\nabla(u_j-u_k)}^2\vp^2\, dX\le C_\varphi \iint  A_j \nabla(u_j-u_k)\cdot \nabla(u_j-u_k)\vp^2\,dX \\
 .
\end{multline*}
By a Caccioppoli-type argument using both the equations $L_ju_j=0$ and $L_ku_k=0$, one can show that for $j$ and $k$ sufficiently large,
\begin{multline*}
    \iint\abs{\nabla(u_j-u_k)}^2\vp^2\, dX\le C\iint\abs{A_j - A_k}^2\abs{\nabla u_k}^2\vp^2\, dX +  C\iint \abs{u_j-u_k}^2\abs{\nabla \vp}^2dX.
\end{multline*}
see for instance the proof of \cite[Lemma 3.23]{DLMsmall} for details. The second term on the right-hand side tends to 0 by the uniform - hence $L^2$ - convergence of $u_j$. The first term also converges to 0 by the uniform convergence \eqref{AjtoAinfty} on $\supp \varphi$ and the fact that $\iint_{\supp \varphi} |\nabla u_k|^2 dX$ is uniformly bounded (thanks to the Caccioppoli inequality and the uniform convergence of $u_k$). The claim \eqref{eq.localW12conv} follows.

From here, we easily deduce that $u_\infty$ is a weak solution to $L_\infty u = 0$ in $\Omega_\infty$. Indeed,  for $\varphi \in C_0^\infty(\Omega_\infty)$ and $j$ large enough so that $\supp \varphi \Subset \Omega_j$, we have
\begin{multline*}
\left |\iint_{\Omega_\infty} A_\infty \nabla  u_\infty \cdot \nabla \varphi \, dX \right|
=  \left |\iint_{\Omega_\infty} A_\infty \nabla  u_\infty \cdot \nabla \varphi  \, dX  - \iint_{\Omega_j} A_j \nabla  u_j \cdot \nabla \varphi  \, dX \right| \\
\leq \left| \iint_{\Omega_\infty} A_\infty \nabla (u_\infty-u_j) \cdot \nabla \varphi \, dX \right|
+ \left| \iint_{\Omega_\infty} (A_\infty- A_j) \nabla u_j \cdot \nabla \varphi \, dX \right| \rightarrow 0 
\end{multline*} 
by \eqref{eq.localW12conv} and \eqref{AjtoAinfty}.

It remains to show that $\A_\infty$ and $|\nabla u_\infty|^2$ are constant. Take $X\in \Omega_\infty$. Since $\Omega_\infty$ is uniform, it is an open connected domain, and thus there exists a connected open set $S_{X}$ with smooth boundaries such that $X,X_0 \in S$ and $\overline{S_{X}} \Subset \Omega_\infty$. For $r>0$ small enough so that $B(X,r) \cup B(X_0,r) \subset S_X$, we have 
\begin{multline*}
\left| \{\A_\infty\}_{B(X,r)} - \{\A_\infty\}_{B(X_0,r)} \right| \\ 
\leq \{|\A_\infty - \A_j|\}_{B(X,r)} + \{|\A_\infty - \A_j|\}_{B(X_0,r)}
+ C_X \iint_{S_X} \Big| \A_j -  \{\A_j \}_{S_X} \Big| dY \\
=: A_1 + A_2 + A_3.
\end{multline*} 
where we write $\{a\}_{E}$ for $\fint_E a\, dY$. The terms $A_1$ and $A_2$ converge to 0, since $\A_j$ converges uniformly to $\A_\infty$ in compact sets of $\Omega_\infty$. Moreover, for $j$ large enough so that $S_X \subset \Omega_{j,j}$,
\[A_3 \leq C_X \iint_{S_X} |\nabla \A_j | dY \leq C_X 2^{j} \to 0\]
by the Poincar\'e inequality on $S_X$ and the assumption \eqref{assumptionblowup2}. So we proved that $\{\A_\infty\}_{B(X,r)} - \{\A_\infty\}_{B(X_0,r)} = 0$ for all $X\in \Omega_\infty$ and all $r$ small enough, hence $\A_\infty$ is essentially constant. The proof for $|\nabla u_\infty|^2$ is similar. We have that
\begin{multline*}
\left| \{|\nabla u_\infty|^2\}_{B(X,r)} - \{|\nabla u_\infty|^2\}_{B(X_0,r)} \right| \\ 
\leq 2 \Big( \{|\nabla (u_\infty-u_j)|^2\}_{B(X,r)} + \{|\nabla (u_\infty-u_j)|^2\}_{B(X_0,r)} \Big)
+ C_X \iint_{S_X} \Big| |\nabla u_j|^2 -  \{|\nabla u_j|^2 \}_{S_X} \Big| dY \\
=: U_1 + U_2 + U_3.
\end{multline*} 
The terms $U_1$ and $U_2$ converge to 0, since we already proved that $u_j \to u_\infty$ in $W^{1,2}_{loc}(\Omega_\infty)$. As for $U_3$, observe that
\[U_3 \leq C_X \left(\iint_{S_X} \frac{|\nabla\br{|\nabla u_j|^2} |^2}{u_j^4} dY\right)^\frac12 u_j^2(X_0) \leq C_X 2^{j} \to 0\]
by the Poincar\'e inequality on $S_X$, the Harnack inequality (that gives  that $u(Y) \approx u(X_0)$ on $S_X$), and - for $j$ large enough - the assumption \eqref{assumptionblowup}. Altogether, we showed that for any $X\in \Omega_\infty$ and any $r$ small enough,
\[\fiint_{B(X,r)}|\nabla u_\infty|^2 dY= \fiint_{B(X_0,r)} |\nabla u_\infty|^2 dY\]
Since $u_\infty$ is a solution of a constant coefficient elliptic operator, $u_\infty \in C^\infty(D_\infty)$ and hence, by taking $r\to 0$, we obtain $|\nabla u_\infty(X)| = |\nabla u_\infty(X_0)|$ as desired. The proposition follows.
\ep

As a direct consequence of Proposition \ref{problowup2} and Corollary \ref{corDEM}, we obtain
\begin{corollary} \label{problowup3}
Let  $\{\Omega_j\}_{j\geq 1}$, $\{L_j\}_{j\geq 1}$, and $\{u_j\}_{j\geq 1}$ be as in Propositions \ref{problowup1} and \ref{problowup2}.

Then, $d$ is an integer, and up to a subsequence and up to rotation and translation, $\partial \Omega_j$ converges to a $d$-dimensional plane. 
\end{corollary}

\subsection{The blow-up argument}

\begin{definition}
Let $\Omega$ be a domain with $d$-Ahlfors regular boundary $\partial \Omega$ (and we call $\sigma$ the Ahlfors regular measure on $\partial \Omega$). We say that a collection $\B \subset \D_{\partial \Omega}$ is Carleson-packing if there is a constant $C>0$ such that for any $Q_0\in \D_{\partial \Omega}$,
\[ \sum_{ Q \subset Q_0 \atop Q \in \B} \sigma(Q) \leq C \sigma(Q_0).\]
The complement of a Carleson-packing collection $\B$ will always be denoted by $\G$, that is, $\G=\D_{\pom}\setminus\B$.
\end{definition}

The definition has the following useful property.

\begin{lemma} \label{lemcoherent}
Let $\Omega$ be a domain with $d$-Ahlfors regular boundary $\partial \Omega$ (and we call $\sigma$ the Ahlfors regular measure on $\partial \Omega$), and let  $\B \subset \D_{\partial \Omega}$ be Carleson packing. Then $\D_{\partial \Omega} = \B \cup \G$ can be further partitioned  as 
\[ \G = \bigcup_i \cS_i,\]
where each $\cS_i$ is semi-coherent - that is 
\begin{enumerate}
\item $\cS_i$ contains a unique maximal element $Q(\cS_i)$, that is $Q \subset Q(\cS_i)$ for all $Q\in \cS_i$,
\item if $Q \subset R \subset Q(\cS_i)$, then $R \in \cS_i$,
\end{enumerate}
and the collection $\{Q(\cS_i)\}_i$ is Carleson-packing.
\end{lemma}

\bp
It is classical. The ``roots'' $Q(\cS_i)$ are the children of elements in $\B$ that are not themselves in $\B$. We construct $\cS_i$ by induction: $\cS_i$ contains the elements of $\D_{\partial \Omega,k}$ that are not in $\B$ and whose ancestor is in $\cS_i$. The fact that $\{Q(\cS_i)\}_i$ is Carleson-packing comes from the fact that they are children of elements of $\B$, which is itself Carleson-packing. 
\ep

\begin{proposition} \label{propCarlpack}
Let $d<n$, and let $\Omega$ be a domain in $\R^n$ such that its boundary $\partial \Omega$ is $d$-Ahlfors regular. Take $\beta >0$ and define $D_\beta$ as in \eqref{defDbeta}. For $Q\in\D_{\pom}$, define $\Omega^{j}(Q)$ as 
\[
\Omega^j(Q):= \{X\in \Omega \cap B_{Q^{(j)}} , \, \delta_{\pom}(X) \geq 2^{-j}\ell(Q)\},
\]
where $Q^{(j)}$ is the ancestor of $Q$ from $j$ generations before. We have the following properties.
\begin{enumerate}
\item If $\B_x$ is Carleson-packing, then 
\[\B_x(j) := \Big\{Q \in \D_{\partial \Omega}, \, Q^{(k)} \in \B_x \text{ for some } k\in \{0,\dots,j\} \,   \Big\}\]
is Carleson-packing.
\item If $L=-\diver [D_{\beta}^{d+1-n} \A \nabla]$ is a uniformly elliptic operator with $\A$ satisfying $D_\beta \nabla \A \in CM_\Omega(M_\A)$, then the collection 
\[\B_\A(j) := \Big\{Q \in \D_{\partial \Omega}, \, \iint_{\Omega^j(Q)} |\nabla \A| dX \geq 2^{-j} \ell(Q)^{n-1} \Big\}\]
is Carleson-packing.
\item If there is a Carleson-packing family $\B_x$ of dyadic cubes in $\partial \Omega$ and a constant  $M_u>0$ such that for any $Q\in \G_x$, there exists a positive weak solution $u = u_{Q}$ to $Lu = 0$ in $B_Q \cap \Omega$ that satisfies $u = 0$ on $\Delta_Q$ and 
\begin{equation} \label{NNuisCM2}
\iint_{B(y,t) \cap \Omega} \dfrac{|\nabla |\nabla u_{Q}|^2|^2}{|u_{Q}|^4}D_\beta^{d+6-n} \, dX \leq (M_u)^2 t^d,
\end{equation}
whenever $B(y,2t)\subset B_Q$, then the set 
\begin{multline*}
\B_u(j) := \B_x(j+2) \cup \Big\{ Q \in \G_x(j+2), \,
\iint_{\Omega^j(Q)} \frac{|\nabla |\nabla u_{R}|^2|^2}{|u_{R}|^4} dX \geq 4^{-j} \ell(Q)^{n-6}\\
\text{holds for $R \in \G_x$ such that $k(R) \leq k(Q)-j-2$ and $Q \in \D_{\partial \Omega}(R)$}\Big\}
\end{multline*}
is Carleson-packing.
\end{enumerate}
\end{proposition}

\bp
The first collection is easy. Since $\sigma(Q) \leq C2^{-kd} \sigma\br{Q^{(k)}}$, we have
\[\sum_{ Q \subset Q_0 \atop Q \in \B_x(j)} \sigma(Q) \lesssim \sum_{k=0}^j 2^{-kd}  \sum_{ Q \subset (Q_0)^{(k)} \atop Q \in \B_x} \sigma(Q) \lesssim \sum_{k=0}^j 2^{-kd} \sigma\br{Q_0^{(k)}} \lesssim C(j+1) \sigma(Q_0).\] 

For the second collection $\B_\A(j)$, observe that for $Q\in \B_\A(j)$, we have
\[ 1 \leq \left(2^j \ell(Q)^{1-n} \iint_{\Omega^j(Q)} |\nabla \A| dX\right)^2 \leq C 2^{j(2+n)}  \ell(Q)^{-d} \iint_{\Omega^j(Q)} |\nabla \A|^2 D_\beta^{d+1-n} \, dX,\]
thanks to \eqref{Db=dist}. Therefore
\[\sigma(Q) \leq C_j \iint_{\Omega^j(Q)} |\nabla \A|^2 D_\beta^{d+1-n} \, dX\]
and then, if $Q^*_0$ is the ancestor of $Q_0$ from $j$ generations ago, 
\begin{multline*}
\sum_{Q \subset Q_0 \atop Q \in \B_\A(j)} \sigma(Q)
\leq C_j \sum_{Q \subset Q_0 \atop Q \in \B_\A(j)}  \iint_{\Omega^j(Q)} |\nabla \A|^2 D_\beta^{d+1-n} \, dX
\leq C \iint_{B_{Q^*_0} \cap \Omega} |\nabla \A|^2 D_\beta^{d+1-n} \, N_j(X) \, dX,
\end{multline*}
where $N_j(X)$ is the number of dyadic cubes such that $X\in \Omega^j(Q)$. Since $X\in \Omega^j(Q)$ implies $\ell(Q) \leq 2^j\delta(X)$ and $|X-x_Q| \leq 2^{j+1}\ell(Q)$, we deduce that $N_j(X)$ is bounded independently of $X$ (it depends on $j$ but we don't care). As such, we obtain
\[
\sum_{Q \subset Q_0 \atop Q \in \B_\A(j)} \sigma(Q)  \leq C_j  \iint_{B_{Q^*_0} \cap \Omega}  |\nabla \A|^2 D_\beta^{d+1-n} \, dX \leq C_j (M_\A)^2 \ell(Q^*_0)^{d} \leq C'_j (M_\A)^2 \sigma(Q_0)
\]
since $D_\beta \nabla \A \in CM_\Omega(M_\A)$. We deduce that $\B_\A(j)$ is Carleson-packing.

We turn to the third collection $\B_u(j)$. We partition $\G_x(j+2)$ into the semi-coherent regimes $\bigcup_i \cS_i$ given by Lemma \ref{lemcoherent}. Since the $\B_x(j+2)$ and the $\{Q(\cS_i)\}_i$ are Carleson-packing, we just need to prove that
\begin{equation} \label{claimBuj}
\sum_{Q \in \cS_i \cap \B_u(j)} \sigma(Q) \lesssim \sigma(Q(\cS_i)).
\end{equation}
Since $Q(\cS_i) \in \G_x(j+2)$, the $j^{th}$ and $(j+2)^{th}$ ancestors of $Q(\cS_i)$ - we call them $Q^*_i$  and $Q^{**}_i$ - lie in $\G_x$, and thus, for any $Q\in \cS_i \cap \B_u(j)$, we have by definition  $Q \in \B_u(j) \setminus \B_x(j+2)$ and that 
\[ \sigma(Q) \approx \ell(Q)^d \leq 4^j \ell(Q)^{6+d-n} \iint_{\Omega^j(Q)} \frac{|\nabla |\nabla u_{Q^{**}_i}|^2|^2}{|u_{Q^{**}_i}|^4} dX \leq C_j \iint_{\Omega^j(Q)} \frac{|\nabla |\nabla u_{Q^{**}_i}|^2|^2}{|u_{Q^{**}_i}|^4} D_\beta^{6+d-n} dX\]
by \eqref{Db=dist}.  As a consequence,
\begin{multline*}
\sum_{Q \in \cS_i \cap \B_u(j)} \sigma(Q) \leq C_j \sum_{Q \in \cS_i \cap \B_u(j)} \iint_{\Omega^j(Q)} \frac{|\nabla |\nabla u_{Q^{**}_i}|^2|^2}{|u_{Q^{**}_i}|^4} D_\beta^{6+d-n} \, dX  \\
\leq C_j  \iint_{B_{Q^*_i} \cap \Omega}  \frac{|\nabla |\nabla u_{Q^{**}_i}|^2|^2}{|u_{Q^{**}_i}|^4} D_\beta^{6+d-n} \, N_j(X) \, dX,
\end{multline*}
where $N_j(X)$ is as before the number of dyadic cubes such that $X\in \Omega^j(Q)$, which is uniformly bounded in $X$. We conclude that
\[\sum_{Q \in \cS_i \cap \B_u(j)} \sigma(Q)  \leq C_j  \iint_{B_{Q^*_i} \cap \Omega}  \frac{|\nabla |\nabla u_{Q^{**}_i}|^2|^2}{|u_{Q^{**}_i}|^4} D_\beta^{6+d-n} \, dX \leq C_j(M_u)^2 \ell(Q^*_i)^{d} \leq C'_j \sigma(Q(\cS_i))\]
by \eqref{NNuisCM2}, since $2B_{Q^*_i} \subset B_{Q^{**}_i}$. The claim \eqref{claimBuj} and then the proposition follows.
\ep

\begin{theorem} \label{Thconverse}
Let $d<n$, and let $\Omega \subset \R^n$ be a uniform domain with a $d$-Ahlfors regular boundary. For $\beta >0$, define $D_\beta$ as in \eqref{defDbeta}, and take a uniformly elliptic operator $L$ in the form $-\div [D_\beta^{d+1-n} \A \nabla]$ with $\A$ that satisfies $D_\beta \nabla \A \in CM_\Omega$. Assume finally that there exists a Carleson-packing collection $\B_x$ and a constant $M_u>0$ such that for any $Q\in \G_x := \D_{\partial \Omega} \setminus \B_x$, there exists a positive weak solution $u = u_{Q}$ to $Lu = 0$ in $B_Q \cap \Omega$ that satisfies $u = 0$ on $\Delta_Q$ and 
\begin{equation} \label{NNuisCM}
\iint_{B(y,t) \cap \Omega} \dfrac{|\nabla |\nabla u_{Q}|^2|^2}{|u_{Q}|^4}D_\beta^{d+6-n} \, dX \leq (M_u)^2 t^d,
\end{equation}
whenever $B(y,2t)\subset B_Q$. Then $d$ is an integer and $\partial \Omega$ is uniformly rectifiable.
\end{theorem}

\bp We shall use the characterization of uniform rectifiability given in Theorem \ref{ThUR}. Let $\epsilon>0$ and define
\[
\B_{ur}(\epsilon):=\set{Q\in\D_{\pom},\, b\beta_\infty(Q) > \epsilon}.
\]
The theorem will be proved if we can show that $\B_{ur}(\epsilon)$ is Carleson-packing, and thanks to Proposition \ref{propCarlpack}, this will be immediate once we establish that there exists $j = j(\epsilon)\in \mathbb N$ such that
\begin{equation} \label{claimBur}
\B_{ur}(\epsilon) \subset \B_\A(j) \cup \B_u(j).
\end{equation}
Without loss of generality, we can pick a random $x_0\in \partial \Omega$ and add to $\B_x$  all the dyadic cubes that contain $x_0$, since it is easy to check that the collection $\{Q\in \D_{\partial \Omega}, \, Q \ni x_0\}$ is Carleson-packing. Under this assumption, any cube $Q\in \G_u(j)$ satisfies $x_0 \notin B_{Q^{(j)}}$, and thus
\begin{equation} \label{diam>2j}
\ell(Q) \leq 2^{-j}\diam \partial \Omega \quad \text{ whenever } Q \in \G_u(j).
\end{equation}

Assume by contradiction that \eqref{claimBur} is false for all $j$. Then there exists a $\epsilon_0>0$, a collection $\Omega_j$ of uniform domains  with $d$-Ahlfors regular boundaries (with the constants in \eqref{defADR}, Definition \ref{defCP} and Definition \ref{defHC} being uniform in $j$), a collection of uniformly elliptic operators $L_j=-\div [D_{\beta,j}^{d+1-n} \A_j \nabla]$ (with elliptic and boundedness constants that are uniform in $j$) that satisfies $\sup_j \|\delta_{\pom_j} \nabla \A_j\|_\infty < \infty$, a collection of cubes $Q_j \in \D_{\partial \Omega_j}$, and a collection of weak solutions $u_j$ to $L_ju_j = 0$ in $B_{Q_j^{(j+2)}} \cap \Omega_j$ that satisfies $u=0$ on $\partial \Omega_j \cap B_{Q_j^{(j+2)}}$, such that 
\begin{equation} \label{contradiction1}
\iint_{\Omega_{j,j}(Q_j)} |\nabla \A_j| dX \leq 2^{-j} \ell(Q_j)^{n-1}
\end{equation}
and
\begin{equation} \label{contradiction2}
\iint_{\Omega_{j,j}(Q_j)} \frac{|\nabla |\nabla u_j|^2|^2}{|u_j|^4} dX \leq 4^{-j} \ell(Q_j)^{n-6}
\end{equation}
but
\begin{equation} \label{contradiction3}
b\beta_\infty(Q_j) > \epsilon.
\end{equation}
Here $Q^{(j+2)}_j$ is the ancestor of $Q_j$ from $j+2$ generations before and 
\[\Omega_{j,j}(Q_j) := \{X\in \Omega_j \cap B_{Q_j^{(j)}} , \, \delta_{\pom_j}(X) \geq 2^{-j}\ell(Q_j)\}.\]

By dilatation and translation invariance, we can assume that $0= x_{Q_j^{(j)}}$ - i.e. $0$ is the center of $B_{Q_j^{(j)}}$ - and $\ell(Q_j) = 1$. Moreover, \eqref{diam>2j} shows that $\diam \partial \Omega_j\geq 2^j$. So the $\Omega_j$, $L_j$, $u_j$ satisfies the assumptions of Corollary \ref{problowup3}, which entails that $\partial \Omega_j$ converges to an affine $d$-plane $P_\infty$, in particular
\[ \lim_{j\to \infty} \left( \sup\{ \dist(y,P_\infty), \, y \in \partial \Omega_j \cap 2B_{Q_j}\} +  \sup\{ \dist(y,\partial \Omega_j), \, y \in P_\infty \cap B(x,r)\} \right) = 0,\]
and thus
\[  \lim_{j\to \infty} b\beta_\infty(Q_j) = 0,\]
which gives a contradiction with \eqref{contradiction3}. The claim \eqref{claimBur} and the theorem follows.
\ep

The implication $(iv)\implies (i)$ in Theorem \ref{mainthm} is a simple consequence of Theorem \ref{Thconverse}.
\medskip

\noindent {\em Proof of Theorem \ref{mainthm} $(iv)\implies(i)$: }
Given $Y\in\om$, we show that the collection
\[
\B_Y:=\set{Q\subset\D_{\pom}, \, Y\in B_{Q}}
\]
is Carleson-packing. Fix any $Q_0\in \D_{\pom}$ and set $k_0:=k(Q_0)$. If $Y\notin 2B_{Q_0}$, then for any dyadic cube $Q\subset Q_0$, $\abs{Y-x_Q}\ge \abs{Y-x_{Q_0}}-\abs{x_{Q_0}-x_Q}\ge \ell(Q_0)\ge2\ell(Q)$, and thus $Q\notin \B_Y$; so we of course have 
\(
\sum_{Q\subset Q_0, Q\in\B_Y}\sigma(Q) = 0 \le \sigma(Q_0)
\).
In the case where $Y\in 2B_{Q_0}$, we observe that if $Q,Q'\in\D_{k}$ and $Y\in B_Q\cap B_{Q'}$, then $\abs{x_Q-x_{Q'}}\le2^{k+1}$. Since the cubes of generation $k$ are disjoint, it entails that the number of the cubes $Q\in\D_k$ for which $Y\in B_Q$ is finite and depends only on the dimension and the Ahlfors-regular constant. So we estimate crudely and obtain
\[
\sum_{Q\subset Q_0 \atop Q\in\B_Y}\sigma(Q)\le\sum_{j=k_0}^\infty\sum_{Q\in\D_j\atop Q\in\B_Y}\sigma(Q)\le C\sum_{j=k_0}^\infty 2^{-jd}\le C2^{-k_0d}\le C\sigma(Q_0).
\]
For $Q\in\G_Y$, the Green function $G^Y$ is a positive solution in $B_Q\cap\om$ and satisfies $G^Y=0$ on $\Delta_Q$. So Theorem \ref{Thconverse} applies and we get that $d$ is an integer and $\pom$ is uniformly rectifiable. \ep

\bibliographystyle{alpha}
\bibliography{reference}
\end{document}